\documentclass[11pt]{amsart} 

\usepackage[all]{xy}
\usepackage{amsthm}
\usepackage{mathrsfs}
\usepackage{amssymb}
\usepackage{amsmath}
\usepackage{array, url}
\usepackage{stmaryrd}
\usepackage{mathtools}
\usepackage{wasysym}
\usepackage{tensor}
\usepackage{mathdots}
\usepackage{amsbsy}

\usepackage[OT2,T1]{fontenc}
\DeclareSymbolFont{cyrletters}{OT2}{wncyr}{m}{n}
\DeclareMathSymbol{\Sha}{\mathalpha}{cyrletters}{"58}

\usepackage{latexsym}
\usepackage{amsxtra}
\usepackage{sagetex}

\usepackage[margin=1in]{geometry}

\let\oldmarginpar\marginpar
\renewcommand\marginpar[1]{\-\oldmarginpar[\raggedleft\footnotesize #1]%
{\raggedright\footnotesize #1}}

\newcounter{firstnumber}[section]
\newcounter{secondnumber}[firstnumber]
\newcounter{thirdnumber}[secondnumber]
\newcounter{fourthnumber}[thirdnumber]
\newcounter{fifthnumber}[fourthnumber]
\newcounter{currentdepth}

\renewcommand{\thefirstnumber}{\arabic{section}.\arabic{firstnumber}}
\renewcommand{\thesecondnumber}{\thefirstnumber.\arabic{secondnumber}}
\renewcommand{\thethirdnumber}{\thesecondnumber.\arabic{thirdnumber}}

\newcommand{\segment}[2]{%
    \setcounter{currentdepth}{1}%
    \def\thesubsection{\thefirstnumber}%
    \refstepcounter{firstnumber}\label{#1}%
    \addtocounter{subsection}{-1}%
    \subsection{#2}}
\newcommand{\ssegment}[2]{%
    \setcounter{currentdepth}{2}%
    \def\thesubsection{\thesecondnumber}%
    \refstepcounter{secondnumber}\label{#1}%
    \addtocounter{subsection}{-1}%
    \subsection{#2}}
\newcommand{\sssegment}[2]{%
    \setcounter{currentdepth}{3}%
    \def\thesubsection{\thethirdnumber}%
    \refstepcounter{thirdnumber}\label{#1}%
    \addtocounter{subsection}{-1}%
    \subsection{#2}}

\newcommand{\set}[2]{\big\{ #1 \; \big| \; #2 \big\} }

\newcommand{\Vect}{\operatorname{\bf{Vect}}}

\newcommand{\from}{\leftarrow}
\newcommand{\xto}{\xrightarrow}
\newcommand{\xfrom}{\xleftarrow}
\newcommand{\surj}{\twoheadrightarrow}

\newcommand{\gr}{\operatorname{gr}}
\newcommand{\Gr}{\operatorname{Gr}}

\newcommand{\Spec}{\operatorname{Spec}}

\newcommand{\Li}{\operatorname{Li}}

\renewcommand{\Im}{\operatorname{Im}}

\newcommand{\Hom}{\operatorname{Hom}}
\newcommand{\Ext}{\operatorname{Ext}}
\newcommand{\Isom}{\operatorname{Isom}}
\newcommand{\End}[1]{\mathrm{End} \, {#1}}

\newcommand{\ad}{\operatorname{ad}}
\newcommand{\Lie}{\operatorname{Lie}}

\newcommand{\m}[1]{\mathrm{#1}}
\newcommand{\fk}[1]{\mathfrak{#1}}

\newcommand{\bb}[1]{\mathbb{#1}}

\newcommand{\la}{\lambda}
\newcommand{\ka}{\kappa}

\newcommand{\ze}{\zeta}
\newcommand{\ga}{\gamma}
\newcommand{\al}{\alpha}

\newcommand{\om}{\omega}

\newcommand{\GL}{\operatorname{GL}}
\newcommand{\gl}{\operatorname{\fk{gl}}}
\newcommand{\Ga}{{\mathbb{G}_a}}
\newcommand{\Gm}{{\mathbb{G}_m}}
\newcommand{\Qp}{{\QQ_p}}
\newcommand{\Zp}{{\ZZ_p}}
\newcommand{\Fp}{{\FF_p}}

\newcommand{\ZZ}{\bb{Z}}
\newcommand{\CC}{\bb{C}}

\newcommand{\nN}{\fk{n}}

\newcommand{\NN}{\bb{N}}
\newcommand{\QQ}{\bb{Q}}

\newcommand{\PP}{\bb{P}}

\newcommand{\Uu}{\mathcal{U}}
\newcommand{\Bb}{\mathcal{B}}
\newcommand{\FF}{\bb{F}}
\renewcommand{\AA}{\bb{A}}

\newcommand{\Oo}{\mathcal{O}}

\newcommand{\SSS}{{S}}

\newcommand{\mM}{\fk{m}}

\newcommand{\inv}{^{-1}}

\newcommand{\areq}{\ar@{=}}
\newcommand{\suphook}{\ar@{^(->}}
\newcommand{\subhook}{\ar@{_(->}}

\newcommand{\inj}{\hookrightarrow}

\newcommand{\iso}{\stackrel{\sim}{\to}}

\newcommand{\abs}[1]{\lvert#1\rvert}

\newcommand{\thrpl}{\PP^1 \setminus \{0,1,\infty\}}

\newcommand{\Mot}{{\bf Mot}}
\newcommand{\Fphi}{\boldsymbol {F\phi} }

\newcommand{\dR}{{\rm {dR}}}

\newcommand{\et}{{\textrm {\'et}}}

\newcommand{\ZSinv}{\ZZ[S\inv]}

\newcommand{\Sbar}{\overline S}
\newcommand{\SSSbar}{\bar \SSS}

\newcommand{\Un}{U^n}
\newcommand{\Unom}{U^{n, \om}}
\newcommand{\UnFphi}{U^{n,F\phi}}

\newcommand{\Sel}{\operatorname{ Sel} }

\newcommand{\bS}{{\bf S}}

\usepackage{etoolbox}
\makeatletter
\patchcmd{\@settitle}{\uppercasenonmath\@title}{}{}{}
\makeatother

\begin{document}

\SelectTips{cm}{11}

\title{Mixed Tate motives and the unit equation}
\thanks{Both authors are supported by the DFG priority programme SPP 1489 \textit{Algorithmic and experimental methods in Algebra, Number Theory and Geometry}.}

\author{Ishai Dan-Cohen}
\author{Stefan Wewers}

\date{\today}

\begin{abstract}
This is the second installment in a sequence of articles devoted to ``explicit Chabauty-Kim theory'' for the thrice punctured line. Its ultimate goal is to construct an algorithmic solution to the unit equation whose halting will be conditional on Goncharov's conjecture about exhaustion of mixed Tate motives by motivic iterated integrals (refined somewhat with respect to ramification), and on Kim's conjecture about the determination of integral points via $p$-adic iterated integrals. In this installment we explain what this means while developing basic tools for the construction of the algorithm. We also work out an elaborate example, which goes beyond the cases that were understood before, and allows us to verify Kim's conjecture in a range of new cases. 


\end{abstract}

\maketitle

%%%%%%%%%%%%%%%%%%%%%%%
\section{Introduction}%%%%%%%%%%%%
%%%%%%%%%%%%%%%%%%%%%%%%%

\segment{rs1}{}
Let $S$ be an open subscheme of $\Spec \ZZ$, and let $X \to S$ be a regular minimal model of a hyperbolic curve over $\QQ$.  The theorems of Siegel and Faltings in this case can be summarized by saying that if the fundamental group of the Betti topological space of $X_\CC$ is nonabelian, then the set $X(S)$ of integral points of $X$ is finite. As soon as we know that $X(S)$ is finite, we're faced with the problem of computing it explicitly. Available in only very special cases, the quest for an \emph{effective Mordell's conjecture} remains one of the holy grails of arithmetic geometry.  

 Starting with his groundbreaking works \cite{kimi, kimii}, Minhyong Kim has developed a new approach to the study of integral points, which seeks to make effective use of the nonabelian nature of the fundamental group to bound, and hopefully compute, the set of integral points. The hope to do so was already inherent in Grothendieck's section conjecture concerning the profinite fundamental group. So far however, the context of the section conjecture hasn't yielded effective methods, even conjecturally. Instead, Kim passes to the \emph{unipotent fundamental group}. Developed in works of Deligne \cite{Deligne89}, Wojtkowiak \cite{Wojtkowiak}, Goncharov \cite{GonMPMTM}, Deligne-Goncharov \cite{DelGon}, Levine \cite{LevineTMFG}, and others, the unipotent fundamental group is actually not one group, but rather a system of groups, together with transcendental \emph{comparison isomorphisms} which reflect its conjecturally motivic origin. Being prounipotent, it lives somewhere between the profinite fundamental group and the abelian world of motives. 
 
 \segment{rs2}{}
  Of particular importance in Kim's original formulation is its $p$-adic \'etale realization. Let $x \in X(S)$ be an integral point, and let $p \in S$ be a prime of good reduction. The theory of the unipotent $p$-adic \'etale fundamental group assigns to $x$ a prounipotent $\Qp$-group $U^\et_x$ equipped with an action of the total Galois group $G_{\QQ}$, and to each $y \in X(S)$ a $G_\QQ$-equivariant $U^\et_x$-torsor $_yP_x^\et$. Kim uses the theory of the unipotent fundamental group to bound the number of points of $X(S)$ as follows. Let $U_x^{n, \et}$ denote the quotient by the $(n+1)^\m{st}$ step in the descending central series, and let $_yP_x^{n, \et}$ denote the associated $U_x^{n,\et}$-torsor. Kim constructs a finite type affine $\Qp$-scheme $\Sel^n_S(X)$ which parametrizes $G_\QQ$-equivariant $U_x^{n,\et}$-torsors obeying suitable local ``selmer conditions'' over $S$, as well as a certain local variant $H^1_f(G_\Qp, U_x^{n,\et})$. The assignment
\[
y \mapsto {_yP_x^{n,\et}}
\]
defines a map $\ka: X(S) \to \Sel_S^n(X)$. Together with its local analog, this map fits into a commuting square
\[
\xymatrix{
X(S) 
\ar[d]_-{\ka} \ar[r] & X(\Zp) \ar[d]^-{\ka^p} \\
\Sel^n_S(X) \ar[r]_-{\m{loc}} & H^1_f(G_\Qp, U_x^{n,\et}).
}
\tag{$*$}
\]
%qqq
We set $X(\Zp)_{S,n}:= (\ka^p)\inv(\Im \m{loc})$. As $n$ varies, these sets form a nested sequence
\[
X(\Zp) \supset X(\Zp)_{S,1} \supset X(\Zp)_{S,2} \supset \cdots \supset X(S).
\]
As soon as $X(\Zp) \neq X(\Zp)_{S,n}$, the sets $X(\Zp)_{S,n}$ become finite, and then the number of points provides a bound for the number of integral points $X(S)$. Moreover, Kim conjectures that this bound is eventually sharp, that is, there is an equality
\[
X(\Zp)_{S,n} = X(S)
\tag{K}
\]
for large $n$ \cite{nabsd}.

\segment{rs3}{}%%%%%%%%
Much of the interest in the construction above and the conjectured equality is due to the belief in the computability of the sets $X(\Zp)_{S,n}$. To make the map $\ka^p$ computationally accessible, Kim turns to the de Rham realization of the unipotent fundamental group. The theory of the de Rham fundamental group applied to $X_\Qp$ assigns to $x$ a prounipotent $\Qp$-group $U_x^\dR$ and to each $y \in X(\Zp)$ a path torsor $_yP_x^\dR$. The path torsor may be interpreted as the space of Tannakian paths for the category of unipotent integrable connections. The equivalence between this category and the category of unipotent overconvergent isocrystals on $X_\Fp$ endows $_yP_x^\dR$ with an action by Frobenius. The ring of functions $\Oo({_yP_x^\dR})$ may alternatively be constructed roughly as the de Rham $0^\m{th}$-cohomology of a certain cosimplicial model $_yX_x^\bullet$ of the path space. This interpretation endows $\Oo({_yP_x^\dR})$ with a lifting of the Hodge filtration to the unipotent setting. Its first step defines a subgroup
\[
F^0U_x^\dR
\]
and sub-torsors
\[
F^0{_yP_x^\dR}.
\]
According to a theorem of Besser \cite{Besser} and Vologodsky \cite{Vologodsky}, each path torsor $_yP_x^\dR$ possesses a unique Frobenius-invariant path $p^\m{Crys}$. Choosing arbitrarily a point $p^\m{Hodge}$ in $F^0{_yP_x^\dR}$, we obtain a map
\[
\al: X(\Zp) \to U_x^{n, \dR} / F^0
\]
which sends $y$ to the unique point $u$ such that
\[
p^\m{Crys} = p^\m{Hodge} \cdot u.
\]
A certain lifting of the Bloch-Kato exponential map to the unipotent level \cite{kimii} places $\al$ and $\ka^p$ in a triangle
\[
\xymatrix{
X(\Zp) \ar[d]_-{\ka^p} \ar[dr]^-\al \\
H^1_f(G_{\Qp}, U_x^{n,\et}) \ar[r]^-\cong_-h & U_x^{n, \dR} / F^0
}
\tag{$*$}
\]
whose commutativity follows from the unipotent $p$-adic Hodge theory of Olsson \cite{OlssonTowards}.

The map $\al$ serves as a $p$-adic analog of the higher Albanese map of Hain \cite{HainHigher}. In coordinates it is given by certain $p$-adic iterated integrals, also known as \emph{Coleman functions} \cite{Besser}. There are fairly well established methods for producing explicit expressions for these on the one hand, and for computing their values to arbitrary $p$-adic precision on the other.

\segment{rs4}{}%%%%%%%%%%
The map $h \circ \m{loc}: \Sel^n_S(X) \to U_x^{n,\dR} / F^0$ is an algebraic map of finite type $\Qp$-schemes. Computing its scheme-theoretic image presents a greater challenge. A range of special cases which may be approached by the methods of nonabelian Galois cohomology are summarized in \cite{nabsd}. The present article is the second in a series, begun in \cite{CKtwo}, whose 
goal, as we now see it, is to approach the case of the thrice punctured line
\[
X = \thrpl
\]
via the theory of motivic iterated integrals. Particular to this case is that $X$ is of mixed Tate type. This means that the cosimplicial model of the path space $_yX_x^\bullet$ gives rise to a simplicial object of the triangulated category of mixed Tate motives. The latter is known to have a natural t-structure. Taking $H^0$ with respect to this t-structure gives us a prounipotent group object $U_x$ (for $x=y$) and torsor-objects $_yP_x$ of the abelian category of mixed Tate motives. The precise construction, which involves truncating and taking a limit, is explained in Deligne-Goncharov \cite{DelGon}. A more elaborate construction (loc. cit.) allows us to replace the integral points $x$, $y$ with integral tangent vectors at the missing points.

\segment{rs6}{}%%%%%%%%%%
The category of mixed Tate motives over $S$ possesses a canonical fiber functor; the associated Tannakian fundamental group is canonically of the form
\[
G(S) = U(S) \rtimes \Gm
\]
with $U(S)$ prounipotent. The unipotent fundamental group $U(X) = U_0(X)$ of $X$ at the tangent vector $\partial/\partial t$ at $0$ is a prounipotent group over $\QQ$ with an action of $G(S)$, and each path torsor $_yP_0$ is a $G(S)$-equivariant $U(X)$-torsor. In this article, we pass to a certain quotient
\[
U(X) \surj U^n = \prod_{i=1}^n \QQ(i) \rtimes \QQ(1)
\]
known as the ``polylogarithmic quotient'', constructed by Deligne \cite{Deligne89}. There is also an analogous quotient of the $p$-adic de Rham fundamental group, which we denote by $U^{n,F\phi}$, a unipotent group object of the category of mixed Tate filtered $\phi$ modules of Chatzistamatiou-\"Unver \cite{ChatUnv}. This Tannakian category also possesses a canonical fiber functor with associated fundamental group
\[
G(\Fphi) = U(\Fphi) \rtimes \Gm.
\]
We thereby obtain a motivic version of diagrams \ref{rs2}($*$) and \ref{rs3}($*$) like so.
\[
\xymatrix{
X(S) \ar[r] \ar[d]_-\ka & X(\Zp) \ar[d] \ar[dr]^\al \\
H^1(G(S), U^n) \ar[r]_-{F\phi} & H^1(G(\Fphi), U^{n,F\phi}) \ar@{=}[r] & U^{n,F\phi}
}
\tag{G}
\]
Our goal, in terms of this diagram, is to construct an algorithm for computing $p$-adic approximations of the sets
\[
\al\inv \bigg( F\phi \Big( H^1 \big( G(\SSS), U^n  \big) \Big)  \bigg)
\]
for $S$, $n$ and $p$ arbitrary. In \textit{Explicit Chabauty-Kim theory}  \cite{CKtwo} we focused on the case $n=2$. Here we lay the foundations that will be needed for the general case with $S$, $n$, and $p$ arbitrary. We then investigate two new special cases, given by $S = \emptyset$, $n$, $p$ arbitrary, and by $S= \{2\}$, $n=4$, $p$ arbitrary.

\segment{}{}%%%%%%%%%%%%%%%%%
We refer the reader to \textit{Explicit Chabauty-Kim theory} for a more thorough discussion of background. As explained in the introduction, and documented in the appendix, the prospects of a direct approach are unclear. The difficulty lies in the relative structurelessness of the nonabelian cohomology variety
\[
H^1 \big( G(\SSS), U^n \big):
\]
although abstractly isomorphic to affine space, it lacks a natural set of coordinates. Even making direct effective use of its known dimension is hard, since the latter depends on Borel's computations of higher $K$-groups \cite{Borel53, Borel77}, which are real-analytic in nature.

\segment{1113.1521}{}%%%%%%%%%%%%%%%%%%%%%%%%%
Our approach is to factor the filtered $\phi$ realization map $F\phi$ through a certain auxiliary vector space. Let $\bar S \subset S$ denote a possibly smaller open subscheme, and let
\[
A(\SSSbar) = \Oo\big(U(\SSSbar)\big)
\]
denote the graded Hopf algebra of functions on $U(\SSSbar)$. Similarly, let
\[
A(\Fphi) = \Oo\big(U(\Fphi)\big)
\]
denote the graded Hopf algebra of functions on $U(\Fphi)$.

\subsection*{Theorem} The filtered $\phi$ realization map $F\phi$ factors through $A(\SSSbar)_1 \times \prod_1^n A(\SSSbar)_i $ as in the following diagram.
\[
\xymatrix{
H^1\big(G(\SSS), U^n \big) \ar[r]^-{F\phi}	\ar[d]_\la
&
H^1 \big( G(\Fphi), \UnFphi \big) \ar@{=}[r] \ar[d]
&
\UnFphi 
\\
A({ \bar S})_1 \times \prod_1^n A({ \bar S})_i \ar[r]_-{F\phi^A}
&
 A(\Fphi)_1 \times \prod_1^n A(\Fphi)_i  \ar[ur]_-{ev_{u\inv}}
}
\tag{L}
\]
Here $F\phi^A$ denotes the map induced by filtered $\phi$ realization on the level of graded Hopf algebras; see Section 5 below for the definitions of the maps $\la$ and $ev_{u\inv}$, as well as a discussion of the maps which form Kim's cutter. 

\segment{}{}%%%%%%%%%%%%%%%%%%%%%%%%
The problem of constructing bases for the vector spaces $A(\SSSbar)_i$ is well known. Candidate elements are provided by Goncharov's theory of motivic iterated integrals (see, for instance, \cite{GonGal}). Goncharov conjectured that in the two extreme cases $\SSSbar = \Spec \ZZ$ and $\SSSbar = \Spec \QQ$, each $A(\SSSbar)_i$ is spanned by iterated integrals on the plane, suitably punctured (c.f. his ICM lectures \cite{GonICM}). The case $\SSSbar = \Spec \ZZ$ is now a theorem due to Francis Brown \cite{Brown}. It is reasonable to hope that $A(\SSSbar)$ is exhausted by iterated integrals in certain intermediate cases as well. One case has been known since before Brown's theorem: the case $\SSSbar = \Spec \ZZ \setminus \{2\}$ follows from Deligne's work \cite{DelMuN}. 

\segment{}{}%%%%%%%%%%%%
After a brief review of free prounipotent groups, our work begins by exploring the possibility that a certain condition on $\SSSbar$, which makes use of the Archimedean absolute value, suffices to ensure exhaustion of $A(\SSSbar)$ by iterated integrals. We offer a modicum of evidence for this possibility coming from level $n=2$ in Proposition \ref{r_15_1}. The \emph{Archimedean condition} is that $\SSSbar$ be of the form
\[
\SSSbar = \Spec \ZZ \setminus \{\mbox{primes}\le m\}
\,.
\]

\subsection*{Proposition}
Suppose $n=1$ or $2$, and suppose $\SSSbar$ obeys the Archimedean condition. Then $A(\SSSbar)_n$ is spanned by iterated integrals on $X$.

\segment{}{}%%%%%%
In a section devoted to $p$-adic periods, we discuss the category of mixed Tate filtered $\phi$ modules over $\Qp$. One of the main goals is the construction of a certain special element $u \in U(\Fphi)$. We also introduce the notion of a \emph{filtered $\phi$ iterated integral}, a filtered $\phi$ analog of Goncharov's motivic iterated integrals. These are compatible with motivic iterated integrals via filtered $\phi$ realization on the one hand, and give rise to Coleman functions by evaluation on the special element $u$ on the other. 

\segment{}{}%%%%%%%%%%%%%%%%%
Next, we turn to the construction of the diagrams \ref{rs6}.(G) and \ref{1113.1521}.(L). In order to construct the map $\la$, we make two preliminary constructions. The first is a proposition in nonabelian cohomology (segments \ref{1107_10}--\ref{1115.1705}); the second is \emph{Deligne's representation} $\rho^D$, a map
\[
U^n \to U^{(n+1)\times(n+1)}
\]
to the group of endomorphisms of
\[
\bigoplus_1^{n+1} \QQ(i)
\]
which are unipotent with respect to the weight filtration (Segment \ref{1115.1707}). According to the proposition, we have
\[
H^1 \big( G(\SSS), U^n \big) = \Hom^\Gm \big( U(\SSS), U^n \big)
\,.
\]
We define $\la$ (Segment \ref{1708}) by sending a $\Gm$-equivariant homomorphism
\[
\rho: G(\SSS) \to U^n
\]
to its composite with a certain set of coordinates on $U^n$ which arise as appropriate matrix entries in Deligne's representation. As the name suggests, the map $ev_{u\inv}$ is given by evaluation at $u\inv\in U(\Fphi)(\Qp)$. The maps $\al$ and $\ka$ need to be projected from the full unipotent fundamental group onto our quotient of choice $U^n$; we give an explicit formula for the projection, and an ensuing formula for $\al$ based on Furusho's computations in Segment \ref{1112.1715}.

\segment{}{}%%%%%%%%%%%%%%%%
Armed with these constructions, we may give a sketch of our future algorithm (\ref{1712}). At each level $n$, it has two main steps. In step one, it searches among the motivic iterated integrals on the affine line, suitably punctured, for bases of $A(\SSSbar)_i$, $i\le n$; in step two, it produces equations for the image of $\la$ in terms of the corresponding coordinates. An important point is that step one will not depend on making effective use of Borel's computations of higher $K$-groups. In this sense, our construction here is less precise than the one made in \cite{CKtwo} where we made effective use of the vanishing of $K_2(\QQ)\otimes \QQ$ to construct a preferred basis of $A(\SSSbar)_2$.

If our program comes to fruition, then we will have at our disposal an algorithm for computing the set $X(\SSS)$. Our algorithm will be usable in practice, but with one caveat --- it will not be guaranteed to halt. The main theorem of \textit{explicit Chabauty-Kim theory for the thrice punctured line}, as we now see it, will have two parts. (1) If the algorithm halts for the input $S$, then the output is equal to the set $X(\SSS)$ of solutions to the $S$-unit equation. (2) If Kim's conjecture holds for $\SSS$, and exhaustion by iterated integrals holds for an open subscheme
\[
\SSSbar \subset \SSS,
\]
then the algorithm halts for the input $S$.

\segment{1126}{}%%%%%%%%%%%%%%%%%%%%%%%
The remainder of the article is devoted to our two examples. The case $\SSS = \Spec \ZZ$, although new, does not depend on the constructions outlined above. Our result does depend, however, on the conjectured nonvanishing
\[
\ze^p(n) \neq 0
\]
for $n$ odd $\ge 3$. An immediate corollary of Proposition \ref{0907b} is the following

\subsection*{Theorem}
Suppose $\SSS = \Spec \ZZ$. Assume $\ze^p(n) \neq 0$ for $n$ odd $\ge 3$. Then the ideal of Coleman functions defining the locus
\[
\al\inv \bigg( F\phi \Big( H^1 \big( G(\SSS), U^\infty  \big) \Big)  \bigg) \subset X(\Zp)
\]
is generated by the functions
\begin{align*}
\log^p z, && \log^p(1-z), && \mbox{and} &&\Li_n^p(z) \mbox{ for } n \mbox{ odd}. 
\end{align*}

\medskip \noindent
This result should be compared with Section 4 of \textit{A nonabelian conjecture} \cite{nabsd}, where Kim's conjecture
\[
\emptyset = \al\inv \bigg( F\phi \Big( H^1 \big( G(\SSS), U^2  \big) \Big)  \bigg)
\]
for $\SSS = \Spec \ZZ$ is shown to hold already in depth $n=2$ for many primes $p$.

\segment{1120.2}{}%%%%%%%%%%%%%%%%%%%%
The case $\SSS = \Spec \ZZ \setminus \{2\}$, $n=4$ ($p$ arbitrary) is substantially more interesting. Segments \ref{1716}--\ref{1112} are devoted to constructing a basis of $A(\SSS)_i$ for each $i \le 4$. A fundamental computational tool for working with the vector spaces $A(\SSS)_i$ is the family of complexes
\begin{align*}
A(m) && = &&
0 \to A_m \to \underset{i,j \ge 1}{\bigoplus_{i+j = m}} A_i\otimes A_j \to \cdots \to A_1^{\otimes m} \to 0
\end{align*}
as well as the formula
\[
H^i(A(m)) = \Ext^i(\QQ(0), \QQ(m))
\]
(see \cite[Segment 3.16]{BGSV}), which holds in any mixed Tate category.
% in the case at hand, we have
%$
%\Ext^i(\QQ(0), \QQ(m)) = 0
%$
% for $i \ge 2$, $\Ext^1(\QQ(0), \QQ(m)) = 0 $ for $m$ even, $\Ext^1(\QQ(0),\QQ(1)) = \QQ^S$, and $\Ext^1(\QQ(0), \QQ(m)) = \QQ\ze(m)$ for $m$ odd $\ge 3$. 
If we restrict attention to classical zeta values and classical polylogarithms, then Deligne's representation may be used to give a very simple formula for the first boundary map, which is just the \emph{reduced coproduct} (that is, the coproduct minus its components in $A_0 \otimes A_m$ and $A_m \otimes A_0$); compare our formula \ref{0913} with Goncharov's Theorem 1.2 of \textit{Galois symmetries} \cite{GonGal}. A second tool at our disposal is the filtered $\phi$-realization morphism
\[
F\phi: A(\SSS)_i \to A(\Fphi)_i
\,.
\]

\segment{}{}%%%%%%%%%%%%%%%%%
Recall that $\Uu U(X)^{dR} = \Uu U(X)^\om$ is canonically of the form $\QQ\langle\langle x, y \rangle\rangle $. Here $x$ corresponds to monodromy around $0$ while $y$ corresponds to monodromy around $1$. For $w$ a word in $x,y$, and $b$ either an integral point or an integral tangent vector at a puncture, we define $\Li^U_w(b)$ to be the motivic iterated integral of $w$ from $0$ to $b$; see Segment \ref{1733} for a precise definition. We set
\[
\log^U b := \Li^U_x b,
\]
\[
\Li^U_i b := - \Li^U_{x^iy} b
\,,
\]
and
\[
\ze^U(i) := \Li^U_i(1)
\,.
\]
The main result here is  as follows (\ref{1112}).

\subsection*{Proposition} The set $\{ (\log^U 2)^4, (\log^U 2)\ze^U(3), \Li^U_4(1/2) \}$ forms a basis of $A(\Spec \ZZ \setminus \{2\})_4$. 

\medskip \noindent
A large part of what makes the case $n=4$ interesting is the role played by $\Li^U_4$ of an $S$-integral point.

\segment{1120.14}{}%%%%%%%%%%%%%%%%%%%%%%
If we fix arbitrary generators of $U(\SSS)$, these endow both $H^1(G(\SSS), U^n)$ and $A(\SSS)_i$ with coordinates; computing the image of
\[
\la: H^1(G(\SSS), U^n) \to A(\SSS)_1 \times \prod_1^4 A(\SSS)_i
\]
in terms of these is a purely formal matter (Segments \ref{r_15_2}--\ref{r_15_3}). More interesting is the problem of rewriting the resulting equations in terms of our concrete polylogarithmic basis (\ref{1736}). In doing so, we make use of the identity
\[
\Li^U_3(1/2) = \widetilde{(7/8)}\ze^U(3) + \frac{1}{6}(\log^U 2)^3
\,,
\]
in which $\widetilde{(7/8)}$ is a rational number $p$-adically close to $7/8$ for several small $p$; see Segment \ref{1737}. The final result consists of a pair of explicit Coleman functions on $X(\Zp)$ (\ref{1738}). 

\subsection*{Theorem} Assume $\ze^p(3) \neq 0$.\footnote{This is known for $p$ regular, and conjectured in general. See Example 2.19(b) of Furusho \cite{FurushoI}.} Suppose $\SSS = \Spec \ZZ \setminus \{2\}$ and $n=4$. Then the ideal of Coleman functions defining the locus
\[
\al\inv \bigg( F\phi \Big( H^1 \big( G(\SSS), U^n  \big) \Big)  \bigg) \subset X(\Zp)
\]
is generated by the functions
\[
F_2(z) = \Li^p_2(z) + \frac{1}{2} (\log^p z)\log^p(1-z) 
\]
and
\[
F_4(z) = 
\Li^p_4(z) 
+ \widetilde{8/7} C^p
(\log^p z)\Li^p_3(z) 
+  \left(
\widetilde{4/21} C^p +\frac{1}{24}
\right) 
(\log^p z)^3\log^p(1-z)
\,,
\]
where
\[
C^p = 
\frac{(\log^p2)^3}{24\ze^p(3)} + \frac{\Li^p_4(1/2)}{(\log^p2)\ze^p(3)}
\,.
\]

\medskip \noindent
The first was already visible in depth $2$; the second is new. As a corollary, we obtain three identities.

\subsection*{Corollary}
Suppose $b = 2$, $1/2$, or $-1$. Then $F_4(b) = 0$.

\segment{}{}%%%%%%%%%%%%%%%%
The case $\SSS = \Spec \ZZ\setminus\{2\}$ provides a new testing ground for Kim's conjecture. Our computations, documented in Section \ref{r_15_4}, provide overwhelming evidence that the conjectured equality (\ref{rs2}.(K)) holds for all primes $p$ in the range $3 \le p \le 29$, after which the computation time becomes too long. The problem of turning our overwhelming numerical evidence into a proof nevertheless requires more work. This will be subsumed under a general study of approximation of roots of $p$-adic functions defined by iterated integrals in the near future.

\segment{}{}%%%%%%%%%%%%%%%%%%
Actually, following arguments given to us by Hidekazu Furusho and by one of the referees, we can prove that $\widetilde{7/8} = 7/8$ in two ways. We relegate this to the appendices in order to emphasize the fact that in working out our future algorithm in the present special case, we don't need to rely on this equality.

\subsection*{Acknowledgements}
We would like to thank Jennifer Balakrishnan for her help in testing the identity involving $7/8$. We would like to thank Giuseppe Ancona, Andr\'e Chatzistamatiou, Marc Levine, and Ismael Soud\`eres for many helpful conversations about the subject of this article. In particular, in the spring of 2013, Marc Levine's large group of postdocs and graduate students chose Francis Brown's work as the topic for that semester's \textit{Motives Seminar} at Essen. The seminar that ensued was much like a master class on mixed Tate motives; the master's influence on our work is plainly visible. The first author would like to thank his Essen host Jochen Heinloth for providing an environment highly conducive to research, as well as Francis Brown for his hospitality during a short visit to the IHES, and for the very helpful conversations which ensued. Finally, we would like to thank both referees for their many helpful comments.

%%%%%%%%%%%%%%%%%%%
\section{Free prounipotent groups}%%%%%%%%
%%%%%%%%%%%%%%%%

\segment{}{}%%%%%%%%%
We make constant use of free prounipotent groups throughout this work. We begin with a brief review of their structure. Let $k$ be a field, and $N$ a set. For any $k$-algebra $R$ we let $R\langle\langle N \rangle\rangle$ denote the ring of noncommutative formal power series in $N$, endowed with the cocommutative Hopf algebra structure with counit given by the natural surjection of rings
\[
\epsilon: R\langle\langle N \rangle\rangle \surj R,
\]
and comultiplication induced by $\Delta(n) = 1\otimes n + n\otimes 1$ for $n \in N$. An element $F\in R\langle\langle N\rangle \rangle$ is \emph{grouplike} if $\epsilon(F) = 1$ and
\[
\Delta(F) = F\otimes F,
\]
and \emph{Lie algebra like} if $\epsilon(F) = 0$ and
\[
\Delta(F) = 1\otimes F + F\otimes 1
\,.
\]

\segment{4_11}{}%%%%%%%%%%
Let $I$ denote the augmentation ideal of $k\langle \langle N \rangle \rangle$, and let $A$ denote the topological dual
\[
k\langle \langle N \rangle \rangle^\lor = \underset \to \lim \big( k \langle \langle N \rangle \rangle / I^n \big) ^\lor.
\]
Let $U$ denote the group
\[
R \mapsto  \big\{ \mbox{grouplike elements in } R\langle \langle N\rangle \rangle \big\}.
\]
Let $\nN$ denote the space of Lie algebra like elements in $k\langle \langle N \rangle \rangle$. Then $U$ is a proalgebraic $k$-group with coordinate ring $\Oo(U) = A$, Lie algebra
\[
\Lie U = \nN,
\]
and completed universal enveloping algebra
\[
\Uu U = k\langle \langle N \rangle \rangle.
\]
We call $U$ the \emph{free prounipotent group on $N$}. It has the implied universal mapping property with respect to (pro)unipotent groups over $k$. In particular, a representation of $U$ on a finite dimensional $k$-vector space $E$ is the same as a nilpotent representation of $\nN$, which is the same as a map of sets
\[
r:N\to \gl E
\]
whose image is contained in a nilpotent subalgebra.

\segment{0907}{Shuffle product}%%%%%%%%%%%%%%%%%%%%%%%%%
%%%%%%
Under the isomorphism of topological vector spaces
\[
A = \Uu U^\lor,
\]
multiplication in $A$ corresponds 
to the shuffle product of linear functionals. We recall what this means. A word $w$ in the alphabet $N$ is said to be a \emph{shuffle} of words $w'$, $w''$ if $w$ may be obtained by interpolating the letters of $w'$ and $w''$ without changing the order of the letters in the individual words. Given $f,g \in A$ and $w$ a word in $N$, we have
\[
(fg)(w) = \underset{\mbox{of } w', w''} {\underset {w \mbox{ a shuffle} }  \sum }f(w')g(w'')
\,.
\]
In particular, if we denote by $f_w$ the element of $A$ dual to the word $w$, then we have 
\[
f_{w'}f_{w''} = \underset{\mbox{of } w', w''}{\sum_{\mbox{shuffles } w}} f_w.
\]
This product is often denoted by the letter $\Sha$. We warn the reader that the isomorphism of vector spaces $\Uu = \Uu^\lor$ induced by the basis, endows $A$ with a second product, the \emph{concatenation product}, which is often denoted by juxtaposition. We will reserve juxtaposition for the shuffle product and avoid using the concatenation product.

\segment{}{}%%%%%%%%%%
If $P$ is a $U$-torsor (necessarily trivial), we let $\Uu P$ denote its universal enveloping module:
\[
\Uu P := \Uu U \times^U P,
\]
a free $\Uu U$-module of rank one.

%%%%%%%%%%%%%%%%%%%%%%%%%%%%%%%%%%%%
\section{Exhaustion of mixed Tate motives by iterated integrals}%%%%%%%
%%%%%%%%%%%%%%%%%%%%%%%%%%%%%%%%%%%%%%

\segment{140316a}{Review of motivic iterated integrals}%%%%%%%%%%%%
Motivic iterated integrals were studied first by Goncharov \cite{GonGal}, and developed by Francis Brown. Although the version introduced by Goncharov is sufficient for our purposes, we begin by reviewing Brown's construction in a way that clarifies, we hope, the relationship between his treatment of multiple zeta values in \cite{Brown}, and his treatment of multiple polylogarithms in more recent works, such as \cite{BrownSingle, BrownICM}. Notationally, we break from the literature in a way that may help the reader remember the relationships between the various kinds of iterated integrals appearing below. Also, for us $_yP_x$ denotes a torsor of paths from $x$ to $y$.  

\ssegment{21a}{}%%%%%%%%
We recall first the notion of a Tannakian matrix entry. If $T$ is a neutral Tannakian category over a field $k$, and $\om$, $\om'$ are two fiber functors, we write
\[
_{\om'} P_\om = \Isom^\otimes(\om, \om')
\]
for the torsor of paths. If $E$ is an object, $v \in E(\om)$ (= the \emph{fiber} of $E$ at $\om$), and $f \in E(\om')^\lor$, then we have the associated matrix entry: $ [E, v,f ]^{\om, \om'} $ denotes the function
\[
_{\om'} P_\om \to \AA^1
\]
given (on a point $p$ with values in an arbitrary $k$ algebra) by
\[
p \mapsto \langle p(v), f \rangle.
\]

\ssegment{140316b}{}%%%%%%%%%%%%%%
Let $\bar S$ denote either an open subscheme of $\Spec \ZZ$ or $\Spec \QQ$. Let $\Mot(\SSSbar)$ denote the category of mixed Tate motives over $\SSSbar$ with $\QQ$-coefficients. $\Mot(\SSSbar)$ possesses a canonical fiber functor
\[
\om: \Mot(\SSSbar) \to \Vect(\QQ)
\]
given by
\[
\om(E) = \bigoplus_i \Hom \big( \QQ(-i), \gr^W_{2i} E \big).
\]

\subsection*{Claim}
The canonical fiber functor $\om$ is canonically isomorphic to the de Rham fiber functor $\dR$.

\begin{proof}
If $E$ is a mixed Tate motive over $\SSSbar$, then $E^\dR$ is a $\QQ$-vector space equipped with Hodge and weight filtrations. From the fact that the associated mixed Hodge structure is mixed Tate, it follows that the map
\[
\bigoplus_i F^iE^\dR \cap W_{2i}E^\dR \to E^\dR
\]
is bijective. Subsequently, the composite
\begin{align*}
F^iE^\dR \cap W_{2i}E^\dR \subset W_{2i}E^\dR \surj \gr_{2i}^W E^\dR
&= \Hom_{\Vect(\QQ)}\big( \QQ(-i)^\dR, \gr^W_{2i} E^\dR  \big) \\
&= \Hom_{\Mot(\SSSbar)}\big( \QQ(-i), \gr^W_{2i} E  \big)
\end{align*}
is bijective, which completes the construction.
\end{proof}

\subsection*{Definition}
We define the \emph{mixed Tate fundamental group of $\SSSbar$} to be the group
\[
G(\SSSbar) = \Gm \ltimes U(\SSSbar)
\]
associated by Tannaka duality to the canonical fiber functor $\om$. We will also add decorations $G_\om(\bar S)$ or $G_\dR(\bar S)$, etc. as needed.

%\ssegment{140316c}{}%%%%%%%%%%%
%Let $S'$ denote a finite subset of $\QQ$, let $w$ denote a word of length $n$ in the elements of $S'$ and let $a,b$ be either $\QQ$-rational points of $\AA^1 \setminus S'$ or tangent vectors to $\AA^1$ at the points of $S'$. We wish to construct an element
%\[
%\int_a^b w|_{U(\Spec \QQ)}
%\] 
%of the $n^\m{th}$ graded piece $A(\Spec \QQ)_n$ of $A(\Spec \QQ)$. This element will always be contained $A(\SSSbar)_n$ for some $\SSSbar \subset \Spec \ZZ$ open; this will hold, for instance, if $\SSSbar$ is contained in the locus of good reduction for $\AA^1 \setminus S'$, $a$, and $b$. We write
%\[
%\int_a^b w|_{U(\SSSbar)}
%\]
%for the iterated integral above, considered as an element of $A(\SSSbar)$.

\ssegment{140316d}{}%%%%%%%%%%%%%
Let $S'$ denote a finite subset of $\QQ$ and let $a$ be either a $\QQ$-rational point of $\AA^1 \setminus S'$ or a nonzero tangent vector to $\AA^1$ at a point of $S'$. We let $U_a = U_a(\AA^1 \setminus S')$ denote the unipotent fundamental group of $\AA^1 \setminus S'$ at $a$. We let ${_bP_a} = {_bP_a}(\AA^1 \setminus S')$ denote the unipotent path torsor. We let $\Uu U_a$ denote the completed universal enveloping algebra. We let $\Uu{_bP_a}$ denote the completed universal enveloping module, a $\Uu U_a$-$\Uu U_b$ bimodule object of $\Mot(\SSSbar)$. 

\ssegment{4_14_b}{}%%%%%%%%%%
The de Rham realizations $U_a^\dR$, $_bP_a^\dR$ were studied extensively for instance by Deligne \cite{Deligne89} (start with segment 15.52). Similar conclusions were reached in a somewhat different way by Goncharov (see in particular Proposition 3.2 of \cite{GonGal}). The path torsor $_bP_a^\dR$ is canonically trivialized by a certain $\QQ$-point $_b1_a \in {_bP_a^\dR(\QQ)}$. In terms of Deligne's Tannakian interpretation of $_bP_a^\dR$, this special path may be constructed as follows: each unipotent flat bundle $(E, \nabla)$ on $\AA^1_\QQ \setminus S'$ extends uniquely to a flat bundle $(\bar E, \bar \nabla)$ on $\PP^1_\QQ$ with log poles along $S' \cup \{\infty\}$; the maps
\[
E(a) \from \bar E(\PP^1) \to E(b)
\]
are bijective because $\bar E$ is trivial, and $_b1_a(E)$ is the composite
\[
E(a) \to E(b)
\,.
\]
Turning to the fundamental group $U_a^\dR$ itself, it is canonically free prounipotent on $S'$. Hence
\[
\Uu U_a^\dR = \QQ \langle\langle S' \rangle \rangle.
\]
The weight filtration is given by
\[
{W_{-2d}}\Uu{U_a^\dR} = \QQ \langle \langle S' \rangle \rangle_{\ge d},
\]
the subspace spanned by words of length not less than $d$. The Hodge filtration is given by
\[
F^{-i}\Uu U_a^\dR = \QQ\langle\langle S'\rangle \rangle _{\le i}.
\]
Finally, the trivialization
\[
\Uu {_bP_a^\dR} = \Uu U_a^\dR
\]
is compatible with Hodge and weight filtrations. In particular, we have
\[
\big(F^0\Uu {_bP_a}(\QQ) \big) \cap {_bP_a}(\QQ) = \{ _b1_a\}
\,.
\]

\ssegment{140316e}{}%%%%%%%%%%%%%%%%
Recall that the space of functions $\Oo({_bP_a^\dR})$ is equal to the topological dual
\[
\big( \Uu{_bP_a^\dR} \big)^\lor = \underset \to \lim \big( \QQ \langle S' \rangle / I^n \big)
\]
of the universal enveloping algebra. We define the iterated integral
\[
\int_a^b w
\]
to be the function on the de Rham path torsor corresponding to the linear functional dual to $w$.

\ssegment{19d}{Definition}%%%%%%%%%%%%
Let
\[
_\dR P_\m{B} (\Spec \QQ) = \Isom^\otimes (\m{B}, \dR)
\]
denote the torsor of Tannakian paths from the Betti to the de Rham realization functor on the category of mixed Tate motives $\Mot(\Spec \QQ)$, and let 
\[
{_\dR A_\m{B}}(\Spec \QQ) = \bigoplus_{i\in\ZZ} {_\dR A_\m{B}}(\Spec \QQ)_i = \Oo \big( {_\dR P_\m{B}} (\Spec \QQ) \big)
\] 
denote the ring of functions, with grading induced by de Rham realization, as in Segment \ref{140316b}.
Given a Betti path $ch$ from $a$ to $b$, we let $o(ch)$ denote the map
\[
{_\dR P_{\m B}}(\Spec \QQ) \to {_bP_a}(\AA^1 \setminus S')^\dR
\] 
given by
\[
p \mapsto p(ch).
\]
By this we mean the following. The universal enveloping bimodule
$
\Uu {_b P_a}(\AA^1 \setminus S')
$ 
is a pro-mixed-Tate-motive. Hence $p$ defines an isomorphism
\[
\Uu {_b P_a}(\AA^1 \setminus S')^B \xto{\cong} \Uu {_b P_a}(\AA^1 \setminus S')^\dR
\]
from its Betti realization to its de Rham realization. It follows from the compatibility of $p$ with tensor products that $p$ takes grouplike elements to grouplike elements. Hence $p(ch)$ belongs to $ {_bP_a}(\AA^1 \setminus S')^\dR$.

Let $w$ be a word in $S'$ of length $n$.
We define
the \emph{motivic iterated integral}
\footnote{We warn the reader that while our definition generalizes Brown's definition of motivic multiple zeta values in \cite{Brown}, it disagrees slightly with his definition of motivic polylogarithms in \cite{BrownICM}: on the one hand, in place of our universal enveloping bi-module, he considers the ring of functions on the path torsor, and on the other hand, his paths go in the opposite direction from ours, from de Rham realization to Betti realization.}
\[
\int_{ch} w \in {_\dR A_\m{B}}(\Spec \QQ)_n
\]
of $w$ over $ch$ to be the composite
\[
{_\dR P_{\m B}}(\Spec \QQ) \xto{o(ch)} {_bP_a}(\AA^1 \setminus S')
\xto {\int_a^b w} \AA^1.
\] 
The same definition is conveniently summarized in the language of Tannakian matrix entries by
\[
\int_{ch} w := \left[ \Uu {_bP_a}(\AA^1 \setminus S'), ch, \int_a^b w \right]^{\m{B}, \dR}.
\]

\ssegment{19e}{}%%%%%%%%%%%%%
The path torsor ${_\dR P_{\m B}}(\Spec \QQ)$ possesses a special complex point $comp \in {_\dR P_{\m B}}(\Spec \QQ)(\CC)$, the \emph{de Rham isomorphism}. The complex period of a motivic iterated integral is the associated classical iterated integral
\[
\int_{ch} c_1 \cdots c_n (comp) = \int _{ch} \frac{dt}{t-c_1} \cdots \frac{dt}{t-c_n}.
\]

\ssegment{19f}{}%%%%%%%%%%
In analogy with Definition 2.7 of Brown \cite{BrownSingle}, we let
\[
{_\dR A_{\m B}}(\Spec \QQ)^+ \subset {_\dR A_{\m B}}(\Spec \QQ)
\]
denote the largest subalgebra such that
\begin{itemize}
\item[(i)] ${_\dR A_{\m B}}(\Spec \QQ)^+$ has weights $\ge 0$, and
\item[(ii)] the action of $G_\dR(\Spec \QQ)$ on ${_\dR A_{\m B}}(\Spec \QQ)$ restricts to ${_\dR A_{\m B}}(\Spec \QQ)^+$.
\end{itemize}
Our analog of Lemma 2.8 of loc. cit. is that ${_\dR A_{\m B}}(\Spec \QQ)^+$ is the subring generated by matrix entries
\[
[ E, v, f]^{\m{B}, \dR}
\]
with $E$ antieffective (i.e. $\Gr_i^W E = 0 $ for $i >0$).

Let
\[
{_\dR P_{\m B}}(\Spec \QQ)^+ = \Spec {_\dR A_{\m B}}(\Spec \QQ)^+.
\]
Since $\Uu{_bP_a}(\AA^1\setminus S')$ is antieffective, the map $o(ch)$ factors through a map 
\[
o(ch)^+: {_\dR P_{\m B}}(\Spec \QQ)^+ \to {_bP_a}(\AA^1 \setminus S')
\]
which is equivariant for the action of the motivic Galois group
$
G_\dR(\Spec \QQ).
$ 

As in the paragraph following Lemma 2.8 of loc. cit., we have
\[
{_\dR A_{\m B}}(\Spec \QQ)^+_0 = \QQ,
\]
so there is a special rational point
\[
_\dR 1_\m{B}^+ \in  {_\dR P_{\m B}}(\Spec \QQ)^+(\QQ)
\]
given by the augmentation ideal $(\Spec \QQ)^+_{>0}$. Since the special de Rham path $_b1_a$ in $\AA^1 \setminus S'$ has a similar description, it follows that
\[
o(ch)^+(_{\dR}1_\m{B}^+ ) = {_b1_a}.
\]

The scheme ${_\dR P_{\m B}}(\Spec \QQ)^+$ together with its action by $G_\dR(\Spec \QQ)$ plays the role of the orbit closure in Brown \cite{Brown}.

\ssegment{21b}{Definition}%%%%%%%%%%%
Let $o({_{\dR}1_\m{B}^+})$ denote the orbit map
\[
U_\dR(\Spec \QQ) \to  {_\dR P_{\m B}}(\Spec \QQ)^+.
\]
Then the composite $o({_{\dR}1_\m{B}^+}) \circ o(ch)^+$ is equal to the orbit map
\[
o({_b1_a}) : U_\dR(\Spec \QQ) \to {_bP_a}(\AA^1 \setminus S').
\]
Let
\[
A_\dR(\Spec \QQ) = \bigoplus_{i \in \NN} A_\dR(\Spec \QQ)_i = \Oo \big( U_\dR(\Spec \QQ) \big).
\]
If $\om$ is a word in the elements of $S'$ of length $n$ as above, we define the \emph{unipotent motivic iterated integral}
\[
\int_a^b \om |_{U(\Spec \QQ)} \in A_\dR(\Spec \QQ)_n
\]
\emph{associated to $a$, $b$, and $\om$} to be the composite
\[
\xymatrix{
U_\dR(\Spec \QQ) 
\suphook[d]_-{o(_\dR1_\m{B}^+)} 
\ar[dr]^-{o({_b1_a})}
\\
{_\dR P_\m{B}}(\Spec \QQ)^+ 
\ar[r]_-{o(ch)^+}  &
{_bP_a}(\AA^1\setminus S')^\dR 
\ar[r]_-{\int_a^b \om} &
\AA^1.
}
\]
In the language of Tannakian matrix entries, we have
\[
\int_a^b \om |_{U(\Spec \QQ)} = \left[ \Uu{_bP_a}(\AA^1 \setminus S'), {_b1_a}, \int_a^b \om \right]^{\dR, \dR}.
\]
Note that the right hand side is a priori a function on $G_\dR(\Spec \QQ)$. However, since $_b1_a$ is in graded degree zero, its restriction to $\Gm$ is identically zero as long as $\om$ is not the empty word $1$.

\ssegment{21b'}{}%%%%%%%%%
Let $L$ denote the Lefschetz period
\[
L = \left[ \QQ(1), \ga, \frac{dz}{z}   \right]^{B, \dR}
\]
where $ \ga \in H_1^B(\CC^*)$ is the standard generator. Then as in Segment 2.8 of Brown \cite{BrownSingle}, we find that inside ${_\dR P_\m{B}}(\Spec \QQ)^+$, ${_\dR P_\m{B}}(\Spec \QQ)$ is the open subscheme defined by $L \neq 0$ and $U_\dR(\Spec \QQ)$ is the closed subscheme defined by $L=0$. The situation is summarized in the following diagram.
\[
\xymatrix
{
U_\dR(\Spec \QQ) \subhook[r]^{o({_\dR 1_\m{B}^+})}  &
{_\dR P_\m{B}}(\Spec \QQ)^+ \ar[dr]^{o(ch)^+} \\
&
{_\dR P_\m{B}}(\Spec \QQ) \suphook[u] \ar[r]_{o(ch)} &
{_bP_a}(\AA^1 \setminus S') \ar[dr]_{\int_a^b\om} \\
\Spec \CC \ar[ur]_{comp} &
&
&
\AA^1
}
\]
The failure of unipotent motivic iterated integrals to have a well-defined complex period is apparent from the diagram. We have chosen our convention for motivic iterated integrals because we feel that the availability of the map $o(ch)^+$ helps clarify the situation somewhat.

%Let 
%\[
%{_\dR A_{\m B}}(\Spec \QQ)^+ = \QQ \oplus {_\dR A_{\m B}}(\Spec \QQ)_{>0}
%\]
%for the grading induced by the $\Gm$-action coming from the action of the fundamental group $G_\dR (\Spec \QQ)$ of the category of mixed Tate motives at the de Rham fiber functor. Let 
%\[
% {_\dR P_{\m B}}(\Spec \QQ)^+ = \Spec {_\dR A_{\m B}}(\Spec \QQ)^+.
%\]
%Since $o(ch)$ is $\Gm$-equivariant, and since 
%\[
%{_bA_a}(\AA^1 \setminus S') := \Oo \big( {_bP_a}(\AA^1 \setminus S') \big)
%\]
%is of the form $\QQ \oplus A_{>0}$, it follows that $o(ch)$ factors through a map 
%\[
%o(ch)^+: {_\dR P_{\m B}}(\Spec \QQ)^+ \to {_bP_a}(\AA^1 \setminus S').
%\]
%This map sends the point $1$ defined by the ideal $A_{>0}$ to the de Rham path $_b1_a$ described above. Let $o(1)$ denote the orbit map
%\[
%o(1): U_\dR(\Spec \QQ) \to  {_\dR P_{\m B}}(\Spec \QQ)^+.
%\] 

\ssegment{1733}{Motivic polylogarithms and special zeta values}%%%%%%%%%%%%%%%%%%%%%%%%%%%%%%%
%Certain special iterated integrals on $\AA^1 \setminus \{0,1\}$ belong to the old tradition of \emph{classical polylogarithms}, lifted, by the construction of segment \ref{21b}, to the motivic setting. 

For $b$ a $\QQ$-rational base point and $w$ a word in $\{0,1\}$ of length $n$, we set 
\[
\Li^U_w(b) := \int_0^b w|_{U(\Spec \QQ)} = \left[ \Uu{_b P_0}(\AA^1 \setminus \{0,1\}), {_b 1_0}, \int_0^b w  \right]^{\dR, \dR}
\,,
\]
an element of $A_\dR(\Spec \QQ)_n$. Here the limit point $0$ stands for the tangent vector $\partial / \partial t$ at $0$. Certain special words $w$ give rise to  \emph{classical polylogarithms}, lifted, by the construction of segment \ref{21b}, to the motivic setting. The traditional notation, with an added superscript $U$, is given by
\[
\log^U(b) :=\Li^U_0(b), 
\]
and for $i \ge 0$,
\[
\Li^U_{i+1}(b) := -\Li^U_{0^i1}(b) 
\,.
\]
The special case with $b$ equal to the tangent vector $-\partial /\partial t$ at $1$ gives rise to the motivic special zeta values
\[
\ze^U(w) := \int_0^1 w|_{U(\Spec \QQ)} 
\]
and
\[
\ze^U(i+1):= \ze^U(0^i1).
\]
We replace the superscripts $U$ with $p$ or $\infty$ to denote the associated $p$-adic or complex polylogarithms or zeta values.

\segment{21c}{}%%%%%%%%%
If $\SSSbar$ is an open subscheme of $\Spec \ZZ$, then there is a surjection
\[
U(\Spec \QQ) \surj U(\SSSbar).
\]
If the unipotent motivic iterated integral $\int_a^b\om|_{U(\Spec \QQ)}$ factors through $\SSSbar$, then we say that $\int_a^b\om|_{U(\Spec \QQ)}$ is \emph{unramified} over $\SSSbar$ and we write
\[
\int_a^b\om|_{U(\SSSbar)}
\]
for the associated function on $U(\SSSbar)$.

\segment{}{Definition} Given ${\SSSbar } \subset \Spec \ZZ$ and $n\in\NN$, we may ask if $A(\SSSbar)_n$ is spanned by elements of the form $\int_a^b w|_{U(\SSSbar)}$. If so, we say that \emph{exhaustion by iterated integrals} holds for $\SSSbar$ at level $n$.

\segment{}{}%%%%%%%%%%%%%%%%%%%%%%
In his ICM lectures \cite{GonICM}, Goncharov states the conjectures that \textit{exhaustion} holds for $\SSSbar = \Spec \ZZ$ and for $\SSSbar = \Spec \QQ$. As mentioned in the introduction, the former is now a theorem due to Francis Brown \cite{Brown}, while the case $\SSSbar = \Spec \ZZ \setminus \{2\}$ follows from Deligne's work \cite{DelMuN}. Moreover, in the latter case we have
\[
S' = \{1, 0, -1\}
\,.
\] 

\segment{1113.1515}{Naive hope}%%%%%%%%%%%%%%%%%%
For $S'$ arbitrary, exhaustion is not likely to hold, but we believe that the following condition on $\SSSbar$ may be sufficient.

\subsection*{Archimedean condition.} \emph{$\SSSbar$ is of the form
$
\SSSbar = \Spec \ZZ \setminus \{\mbox{primes }\le m\}
$
for some $m \in \NN$.}

\medskip \noindent
Moreover, we believe it may be sufficient to restrict attention to the set
\[
S' = \{0,1, \mbox{primes } \le m\}
\]
and to consider only iterated integrals of the form
\[
\int_a^b w |_{U(\SSSbar)}
\]
with $a,b$ nonzero tangent vectors to $\AA^1$ at the points of  $S'$, and $w$ a word in $S'$.

\segment{r_15_1}{}%%%%%%%%%%%%%%%%%
In this direction, we have the following proposition, which follows easily from our work \cite{CKtwo}.

\subsection*{Proposition} Suppose $n =1$ or $2$, and suppose $\SSSbar$ obeys the Archimedean condition. Then exhaustion holds for $\SSSbar$ at level $n$. Moreover, we may set
\[
S' = \{0,1\}
\]
and consider only the motivic logarithms and dilogarithms
$\log^U(b)$, and $\Li^U_2(b)$ 
for $b \in X(\SSSbar)$, as well as products of pairs of logarithms.

\begin{proof}
We know that $A(\SSSbar)_1$ is spanned by $(\log^U q')$'s, for $q'\in \Sbar$. This establishes the first assertion for $n=1$, with no assumption on $\SSSbar$. However, the requirement
\[
S' = \{0,1\}
\]
of the \textit{Moreover} clause places a harsh condition on the integral points which may intervene. Special for $\SSSbar$ obeying the Archimedean condition is that $q'-1$ is a product of powers of primes $\in \Sbar$; hence $q'$ is an $\SSSbar$-valued point of
\[
X= \AA^1 \setminus \{0,1\}
\,.
\] 

Tate's computation of $K_2(\QQ)$ tells us that
\[
A(\SSSbar)_1 \otimes A(\SSSbar)_1 = \ZZ[\Sbar\inv]^*_\QQ \otimes \ZZ[\Sbar\inv]^*_\QQ
\]
is spanned by elements of the forms
\[
u\otimes v + v\otimes u 
\]
and
\[
t\otimes(1-t)
\,.
\]
Here, again, we use the Archimedian condition on $\Sbar$. Now we have an isomorphism
\[
A_2 \iso A_1 \otimes A_1
\]
which sends
\[
uv \mapsto u\otimes v + v\otimes u 
\]
and
\[
-\Li^U_2(t) \mapsto t\otimes(1-t)
\,.
\]
So our generators of $A_1 \otimes A_1$ lift to (products of) motivic iterated integrals on the punctured line. 
\end{proof}

%%%%%%%%%%%%%%%%%%%%%%%%
\section{$p$-adic periods}%%%%%%%%%
%%%%%%%%%%%%%%%%%%%%%%%%%%%

\segment{4_9}{}%%%%%%%%%%%
A mixed Tate filtered $\phi$-module, in its natural context, is a mixed Tate object of the category of weakly admissible filtered $\phi$ modules. Chatzistamatiou-\"Unver \cite{ChatUnv} explain how this definition can be greatly simplified. If we restrict attention to mixed Tate filtered $\phi$-modules over $\Qp$, we can give an even simpler definition, which, moreover, expresses the symmetry between the Hodge and Frobenius structures in this case. 

\bigskip\noindent
\textbf{Definition.}
A mixed Tate filtered $\phi$-module over $\Qp$ is a vector space $E$ over $\Qp$ equipped with an increasing filtration $W$ indexed by $2\ZZ$ called the \emph{weight filtration}, plus two decreasing filtrations denoted $F$ and $F'$ both indexed by $\ZZ$. $F$ is called the \emph{Hodge filtration}; $F'$ is the \emph{filtration induced by Frobenius}. These are required to be opposite to the weight filtration, meaning that the composite maps
\begin{align*}
F^iE \inj E \surj E/W_{2(i-1)}E && F'^iE \inj E \surj E/W_{2(i-1)}E
\end{align*}
are bijective for all $i$. We set
\begin{align*}
E_i:= F^iE \cap W_{2i}E && \mbox{and} && E'_i:= F'^iE \cap W_{2i}E.
\end{align*}

\segment{}{}%%%%%%%%%
Recall that a mixed Tate category over a field $k$ is a Tannakian category over $k$ equipped with a distinguished object $k(1)$ of rank one, such that each simple object is isomorphic to $k(i): = k(1)^{\otimes i}$ for a unique $i \in \ZZ$, and such that
\[
\Ext^1(k(0), k(i)) = 0
\]
if $i \le 0$. Mixed Tate filtered $\phi$ modules over $\Qp$ form a mixed Tate category over $\Qp$, which we denote by $\Fphi(\Qp)$. Morphisms are simply morphisms of vector spaces respecting all three filtrations. The unit object $\Qp(0)$ is given by $\Qp$ in degree $0$ for each filtration. The tensor product is defined by taking tensor product filtrations in the usual way. The Tate object $\Qp(1)$ is $\Qp$ placed in degree $-1$ for the Hodge filtration and for the filtration induced by Frobenius, and in degree $-2$ for the weight filtration. As with any mixed Tate category, there is a canonical fiber functor
\[
\om^\m{can}_{\Fphi}: \Fphi(\Qp) \to \Vect(\Qp)
\]
given by
\[
\om^\m{can}_{\Fphi}(E) = \bigoplus_i \Hom(\Qp(-i), \gr^W_{2i}E).
\]
The associated fundamental group comes equipped with a semidirect product decomposition
\[
G(\Fphi) = \Gm \ltimes U(\Fphi)
\]
with $U(\Fphi)$ prounipotent.

\segment{4_14_e}{}%%%%%%%%%%%%%%
There is a second naturally occurring fiber functor $\om^\m{underlying}$ on $\Fphi(\Qp)$ which sends a mixed Tate filtered $\phi$ module to the underlying vector space. There are also two natural Tannakian paths
\[
\xymatrix{
\om^\m{can}_{\Fphi} \ar@/^8pt/[r]^-p \ar@/_8pt/[r]_-{p'} & \om^\m{underlying}
}
\]
connecting the two natural fiber functors. These are given by taking direct sums of the natural bijections
\[
E_i \xto {p} \Hom(\Qp(-i), \gr^W_{2i}E) \xfrom {p'} E'_i.
\]
Let $u = p \circ (p')\inv$. More explicitly, for each object $E$, $u(E)$ is the unique element of $\GL(E)$, unipotent with respect to the wight filtration, such that 
\[
u(F') = F.
\]
Since $u(\Qp(1)) = 1$, $u$ belongs to the unipotent radical $U(\Fphi)$. 

\segment{}{}%%%%%%%%%%
The completed universal enveloping algebra $\Uu U(\Fphi)$ inherits a grading from the $\Gm$-action. Let $v_i$ denote the component of $\exp u$ in graded degree $i$, $i \in \ZZ_{\le -1}$. Then it is straightforward to check that $U(\Fphi)$ is free prounipotent on the $v_i$'s (\ref{4_11}).

\segment{12_12_a}{}%%%%%%%%%
A realization functor from mixed Tate motives over $\QQ$ which are unramified at $p$ into $\Fphi$ is constructed in Chatzistamatiou-\"Unver \cite{ChatUnv}. We give a brief overview of their construction and its context. The triangulated category of motives over $\Spec \Qp$ has a $p$-adic \'etale realization functor to the derived category of semistable $p$-adic representations of $G_\Qp$ for every prime $p$ \cite{HuberRealVoev}. A realization functor to the derived category of weakly admissible filtered $\phi$-$N$-modules compatible with $p$-adic de Rham cohomology could in principle be constructed by doing $p$-adic Hodge theory on the triangulated level. This approach is closely related to the work \cite{OlssonTowards} of Olsson, and is currently being developed by Deglise and Niziol.

If instead of working over $\Qp$, we work over $\QQ$ (or any number field), and moreover, if we restrict attention to mixed Tate motives, then the triangulated category possesses a canonical t-structure, giving rise to an Abelian category of mixed Tate motives. The realization functor to semistable representations is compatible with the t-structures, so we obtain a realization functor from the abelian category of mixed Tate motives to the abelian category of mixed Tate semistable representations. Now we may apply the semi-stable Dieudonn\'e functor of Abelian $p$-adic Hodge theory to obtain a realization functor to the abelian category of mixed Tate filtered $\phi$-$N$-modules. Chatzistamatiou-\"Unver prove that the realization of a motive which is unramified at $p$ has trivial monodromy, hence belongs to the subcategory of mixed Tate filtered $\phi$-modules.

Finally, in this apropos, let us mention that a direct approach to constructing a realization functor from the triangulated category of unramified mixed motives over $\Qp$ to the derived category of weakly admissible filtered $\phi$ modules is currently under development in work of Brad Drew.

\segment{4_13_a}{}%%%%%%
Let $S'$ denote a finite set of disjoint
%\mar{I'm actually slightly confused right now about what would be a sufficient condition on $S'$, $a$ and $b$ to ensure that $_bP_a$ be unramified at $p$.}
 sections of $\AA^1$ over $\Spec \Zp$. By a $\Zp$-integral base point of $\AA^1 \setminus S'$, we mean either a $\Zp$-point, or a nowhere vanishing tangent vector to a point of $S'$. Consider $\Zp$-integral base points $a$, $b$. Assume $S'$, $a$, $b$ are disjoint. The theory of the de Rham fundamental group of Deligne \cite{Deligne89} and Wojtkowiak \cite{Wojtkowiak} has a $p$-adic variant due to Olsson \cite{OlssonTowards, OlssonBar} which takes values in the category of weakly admissible filtered $\phi$ modules. In the case at hand, the fundamental group $U_a^{F\phi} = U_a^{F\phi}(\AA^1\setminus S')$ is a unipotent group object in $\Fphi(\Qp)$ and the path torsor $_bP_a^{F\phi}$ is a torsor object in $\Fphi(\Qp)$.

\segment{4_14_a}{}%%%
Suppose briefly that the situation of segment \ref{4_13_a} comes from a global situation, given by a finite subset $S' \subset \QQ$ with good reduction over $\SSS \subset \Spec \ZZ$, and by $\SSS$-integral base points $a$, $b$ of $\AA^1 \setminus S'$. Then it follows from Olsson's theory \cite{OlssonTowards} that the filtered $\phi$ path torsor $_bP_a^{F\phi}$ is equal to the filtered $\phi$ realization of the unipotent path torsor $_bP_a$ in $\Mot(\SSS)$.

\segment{4_14_c}{}%%%%%%%%
After forgetting (the filtration induced by) Frobenius, $U_a^{F\phi}$ and $_bP_a^{F\phi}$ are the usual de Rham fundamental group and path torsor. In particular, the comments of segment \ref{4_14_b} apply with $\Qp$ in place of $\QQ$.

\segment{4_14_d}{Filtered $\phi$ iterated integrals}%%%%%
As in the motivic setting, we define the iterated integral
\[
\int_a^b w
\]
to be the function
\[
_bP_a^{F\phi } \to \AA^1_\Qp
\]
induced by the linear functional
\[
\Qp \langle \langle S' \rangle \rangle \to \Qp 
\]
dual to $w$, and the associated \emph{filtered $\phi$ iterated integral}
\[
\int_a^b w|_{U(\Fphi)}
\]
to be the composite of the iterated integral above with the orbit map
\[
o({_b1_a}) : U(\Fphi) \to {_bP_a^{F\phi}}.
\]
In terms of Tannakian matrix entries, we have
\[
\int_a^b w|_{U(\Fphi)} = \left[ \Uu{_bP^{F\phi}_a}(\AA^1_\Qp \setminus S'), {_b1_a}, \int_a^b \om  \right]^{\om^\m{underlying}, \om^\m{underlying}}.
\]
(Note that here ${_bP^{F\phi}_a}(\AA^1_\Qp \setminus S')$ need not be the filtered $\phi$ realization of a global motivic path torsor, as the points $a$, $b$, $S'$ need not be rational.)

\segment{4_14_f}{Compatibility with motivic iterated integrals}%%%%%%%
Motivic and filtered $\phi$ iterated integrals are compatible in the following sense. If our local situation comes from a global situation as in paragraph \ref{4_14_a}, then the diagram
\[
\xymatrix{
U(\SSS)_\Qp \ar[d]_-{o({_b1_a})}  &
U(\Fphi) \ar[d]^-{o({_b1_a})} \ar[l] \\
\big( bP_a^\dR(\AA^1_\QQ \setminus S')\big)_\Qp \ar@{=}[r] \ar[dr]_-{(\int_a^b w)_{\Qp}} &
_bP_a^{F\phi}(\AA^1_\Qp \setminus S') \ar[d]^-{(\int_a^b w)} \\
 &
\AA^1_\Qp
}
\]
commutes. In particular, the similarity between the notation of paragraph \ref{4_14_d} and paragraphs \ref{140316e}, \ref{21b} poses no danger.

\segment{4_14_g}{Compatibility with $p$-adic iterated integrals}%%%%%%%%%%
Under the orbit map
\[
o({_b1_a}): U(\Fphi) \to {_bP_a^{F\phi}},
\]
the inverse special element $u\inv \in U(\Fphi)(\Qp)$ (\ref{4_14_e}) maps to the unique Frobenius invariant path. Indeed, we have
\begin{align*}
\big(F'^0\Uu U(\Fphi) \big) \cap U(\Fphi) 
	&= u\inv \big(F^0\Uu U(\Fphi) \big) \cap U(\Fphi) \\
	&= \big\{ u\inv({_b1_a}) \big\}.
\end{align*}
It follows from Besser's interpretation of Coleman integration \cite{Besser} that the function
\[
b \mapsto \int_a^b w|_{U(\Fphi)}(u\inv)
\]
is a Coleman function on $(\AA^1\setminus S')(\Zp)$.

\segment{r29a}{Filtered $\phi$ polylogarithms}%%%%%%%%%
In particular, if we set
\[
\Li^{F\phi}_w(b) := \int_0^b w|_{U(\Fphi)},
\]
then
\[
\Li^{F\phi}_w(b) (u\inv) = \Li^p_w(b)
\]
is the $p$-adic polylogarithmic value of Furusho \cite{FurushoI}.

%\mar{I'm confused about the dependence on branch of log.}

%%%%%%%%%%%%%%%%%%%%%%%%%%%%%%%%
\section{Blueprints for an algorithm}%%%%%%%%%%%%%%
%%%%%%%%%%%%%%%%%%%%%%%%%%%%%%%%

\segment{1120.1}{}%%%%%%%%%%%%%%%%%%%%%%%%%
Let $S$ denote an open subscheme of $\Spec \ZZ$, and let $\bar S \subset S$ denote a possibly smaller open subscheme. Let $G(\SSS) = \Gm \ltimes U(\SSS)$ denote the fundamental group of the category of mixed Tate motives over $\SSS$ at the canonical fiber functor (or, equivalently, at the de Rham fiber functor). Let $U^n$ denote the group object
\[
U^n := \QQ(1) \ltimes \prod_1^n \QQ(i).
\]
This means that the de Rham realization of $U^n$ is a unipotent group equipped with an action of $G(\SSS)$ in which $U(\SSS)$ acts trivially, and $\Gm$ acts with weight $-i$ on the component $\QQ(i)$ (which is half the motivic weight).

Let $A(\SSS) = \Oo(U(\SSS))$. Fixing a prime $p \in \SSS$, we let $G(\Fphi) = \Gm \ltimes U(\Fphi)$ denote the fundamental group of the category of mixed Tate filtered $\phi$ modules over $\Qp$ at the canonical fiber functor, and let $A(\Fphi) := \Oo(U(\Fphi))$. Let $U^{(n+1) \times (n+1)}$ denote the group of endomorphisms of
\[
V^n := \bigoplus_1^{n+1} \QQ(i)
\]
which are unipotent with respect to the weight filtration. We then have the two diagrams discussed in the introduction; glued together, they look like so.
\[
\xymatrix{
X(\SSS)\ar[d]_\ka \ar[r] & X(\Zp)\ar[dr]^\al \ar[d]_{\ka^{F\phi}}
\\
H^1(G(\SSS), U^n) \ar[d]_\la \ar[r]^{F\phi^H} 
& 
H^1(G(\Fphi), \UnFphi) \ar[d]_{\la^{F\phi}} \ar[r]_-{ev_{u\inv}} 
& 
\UnFphi \suphook[r]_-{\rho^D} 
& 
U^{(n+1)\times(n+1) , F\phi}
\\
A(\SSSbar)_1 \times \prod_1^n A(\SSSbar)_i \ar[r]_-{F\phi^A}
&
A(\Fphi)_1 \times \prod_1^n A(\Fphi)_i \ar[ur]_-{ev_{u\inv}}
}
\]
Our first goal for this section is to discuss each of the arrows appearing in this diagram. To construct the unipotent Kummer maps $\ka$ and $\ka^{F\phi}$ (in Segment \ref{1710} below), we will have to recall the appearance of $U^n$ as a quotient of the unipotent fundamental group $U(X)$ of $X = \thrpl$. If we wish to emphasize de Rham realizations, then the situation may be described in terms of two concrete proalgebraic groups: the Tannakian fundamental group $G(\SSS) = G_\dR(\SSS)$ of the category of mixed Tate motives over $\SSS$ at the de Rham fiber functor, and the Tannakian fundamental group $U(X)^\dR$ of the category of unipotent vector bundles with integrable connection on $X$ at the tangent vector $\partial/\partial t$ at $0$, plus an abstract action of $G_\dR(\SSS)$ on $U(X)^\dR$ which reflects the motivic origin of $U(X)^\dR$. However, for the bottom portion of the diagram, with which we begin, this extra information is not needed.

\segment{1107_10}{A proposition in nonabelian group cohomology}%%%%%%%%%%%%%%%%%%%%%%

A certain proposition in nonabelian group cohomology will yield the construction of $\la$, and with it the proof of Theorem \ref{1113.1521}, as an immediate corollary.
Let $U$, $U'$ be prounipotent groups over a field $k$, which are equipped with a $\Gm$-action such that the associated Lie algebras $\nN$, $\nN'$ are graded in purely negative degrees. Let $G = \Gm \ltimes U$, $G' = \Gm \ltimes U'$. We consider the action of $G$ on $U'$ through the projection
\[
\chi: G \surj \Gm
\,.
\]

\subsection*{Proposition}
Each equivalence class of cocycles $[c]$ in $H^1(G, U')$ contains a unique representative $c_0$ such that
\[
c_0(g) = 1 \mbox{ for each } g \in \Gm;
\]
its restriction to $U$ is a $\Gm$-equivariant homomorphism
\[
c_0|_U: U \to U'.
\]
The map
\[
[c] \mapsto c_0|_{\Gm}
\]
defines a bijection
\[
H^1(G, U') = \Hom^{\Gm}(U, U').
\]

\ssegment{1107}{}%%%%%%%%%%%%%%%%%%%%
For the proof, it is convenient to start with an intermediate construction --- a bijection
\[
H^1(G, U') = \Hom^\chi(G, G')/U'
\,,
\]
which we denote by $c \mapsto \rho_c$. By $\Hom^\chi$, we mean homomorphisms such that the triangle
\[
\xymatrix{
G \ar@{.>}[rr] \ar@{->>}[dr] && G' \ar@{->>}[dl] \\
& \Gm
}
\]
commutes; $U'$ acts by conjugation. Recall that $Z^1(G, U') $ is the set of morphisms of schemes
\[
c:G \to U'
\]
satisfying the cocycle condition
\[
c(g_1g_2) = c(g_1) \big( g_1c(g_2) \big).
\]
We define $\rho_c$ by
\[
\rho_c(g) := c(g)\chi(g).
\]
It is straightforward to check that this defines a $U'$-equivariant bijection
\[
Z^1(G, U') = \Hom^\chi(G,G'),
\]
and hence a bijection as stated.

\ssegment{1108}{}%%%%%%%%%%%%%%%%%%%
We claim that the set of splittings of the projection $G' \surj \Gm$ forms a torsor under $U'$. The latter acts by conjugation. To check transitivity, we apply Segment \ref{1107} with $U = 1$ trivial, to obtain
\[
\{\mbox{splittings}\}/U' = \Hom^\chi(\Gm, G')/U' = H^1(\Gm, U') = *,
\]
the trivial pointed set. Having checked transitivity, we may check freeness against the given splitting $\Gm \subset G'$. Suppose $u'$ is a point of its stabilizer (with values in a $k$-algebra $k'$). Then for all $t\in \Gm$($k''$, a $k'$-algebra), we have
\[
u't{u'}\inv= t,
\]
so 
\[
u' = tu't\inv,
\]
so $u' \in {U'}^\Gm(k'')$, which, according to our condition on the grading, is trivial.

\ssegment{1115.1705}{}%%%%%%%%%%%%%%%%%%%
It follows from Segment \ref{1108} that every orbit for the action of $U'$ on $\Hom^\chi(G,G')$ contains a unique element $\rho$ for which the reverse triangle
\[
\xymatrix{
G \ar@{.>}[rr] && G' \\
& \Gm \ar[ul] \ar[ur]
}
\]
commutes. Using Segment \ref{1107} again, we have
\begin{align*}
H^1(G,U') &= \Hom^\chi(G, G')/U'
\\
	&= \{\mbox{homomorphisms making both triangles commute}\}
\\
	&= \Hom^\Gm(U, U'),
\end{align*}
which concludes the proof of Proposition \ref{1107_10}.

\segment{1115.1707}{Deligne's representation}%%%%%%%%%%%%%%%%%%%
The subrepresentation of the adjoint representation
\[
U^{n+1} \to \GL(V^n)
\]
 given by the Lie ideal
\[
V^n \subset \nN^{n+1} = \Lie U^{n+1}
\]
in which the factor $\QQ(1)$ is eliminated from the product, factors through
\[
U^{n+1} \surj U^n
\,.
\]
We refer to the resulting map
\[
\rho^D : U^n \to U^{(n+1) \times (n+1)}
\]
as \emph{Deligne's representation}. In terms of the natural generators $x$, $y$ of $\nN^{n+1}$, a basis of $V^n$ is given by
\[
\{ x, (\ad x)y, \dots, (\ad x)^ny \}.
\]
For $n\ge2$, $\rho^D$ is faithful. The diagonal maps $\al$ and $ev_{u\inv}$ as we construct them, will actually land in the image $\rho^D(\UnFphi)$.

\segment{}{Proof of Theorem \ref{1113.1521}}%%%%%%%%%%%%%%%%%%
All that remains to be done is to make two simple constructions.

\ssegment{1708}{Construction of $\la$, $\la^{F\phi}$}%%%%%%%%%%%%%%%%%%
In constructing the maps $\la$, $\la^{F\phi}$, as well as the two maps $ev_{u\inv}$, we make constant use of the isomorphism of Proposition \ref{1107_10}, as well as of Deligne's representation. Applied with $U = U(\SSS)$ and $U' = \Unom $ the canonical realization of $\Un$, the proposition reads
\[
H^1(G(\SSS),\Un) = \Hom^\Gm(U(\SSS), \Unom).
\]
Let $\pi_{i,j}$ denote the projection 
\[
\Un \to \Hom(\QQ(-j), \QQ(-i)) = \QQ(j-i).
\]
For $\rho \in \Hom^\Gm(U(\SSS), \Unom)$ we define
\[
\la(\rho) = \big(
(\pi_{-2,-3}) \circ \rho^D \circ \rho,  (\pi_{-1,-2}) \circ \rho^D \circ \rho, (\pi_{-1,-3}) \circ \rho^D \circ \rho, \dots, (\pi_{-1,-n-1}) \circ \rho^D \circ \rho
\big). 
\]
The filtered $\phi$ variant, $\la^{F\phi}$ is defined similarly.

\ssegment{4_15-10}{Remark}
To clarify the definition, we note that the composite
\[
\xymatrix{
\Un \ar[r]^-{\rho^D} \ar@{=}[dr] & U^{(n+1)\times(n+1),F\phi} \ar[d]^{(\pi_{-2,-3}, \pi_{-1,-2}, \pi_{-1,-3}, \dots, \pi_{-1, -n-1}   )}  \\
	& \AA^{n+1}
}
\]
is an isomorphism; moreover, each $\pi_{i,j}\circ \rho^D$ is homogeneous of degree $j-i$. So $\la$, in other words, simply maps $\rho$ to its components under a set of homogeneous coordinates on $U^n$.

\ssegment{1709}{Construction of $ev_{u\inv}$}%%%%%%%%%%%%%%%
As the symbols suggest, both maps $ev_{u\inv}$ are given by evaluation on the inverse special element $u\inv \in U(\Fphi)(\Qp)$ of segment \ref{4_14_e}. Applied to the product
$
A(\Fphi) \times \prod_1^n A(\Fphi)_i,
$
this takes us to affine $n+1$-space; we may then reverse the coordinates of Remark \ref{4_15-10} to obtain a map to $\UnFphi$ as shown
\footnote{Deligne's representation allows us to work with the endomorphism algebra $\End V^n$ where otherwise we would have to work in the enveloping algebra $\Uu \UnFphi$. This will simplify some computations below. With the somewhat awkward formulation here, we pay a small price.}
.

\bigskip \noindent
Since commutativity of the lower portion of the diagram is clear, this completes the proof of Theorem \ref{1113.1521}.

\segment{1710}{Review of $\ka$}%%%%%%%%%%%%
We turn our attention to the upper portion of the diagram. Let $U(X)= U_0(X)$ denote the unipotent fundamental group of $X = \thrpl$ at the usual tangent vector at $0$. Recall the depth filtration $D$ of Deligne-Goncharov \cite{DelGon}, Segment 6.1: we consider the map
\[
U(X) \to \Ga
\]
induced by the natural inclusion $X \subset \Gm$, we let $D^1$ denote its kernel, and $D^n$ the descending central series of $D^1$. Deligne uses the representation $\rho^D$ to establish an isomorphism, necessarily canonical,
\[
U(X) / D^2 = U^\infty := \QQ(1) \ltimes \prod_1^\infty \QQ(i)
\,.
\]
Given $z \in X(\SSS)$, we have the path torsor object $_z P_0$, constructed, for instance, in Deligne-Goncharov \cite{DelGon}. If we push forward to $\Un$ and apply the canonical fiber functor, we obtain a $G(\SSS)$-equivariant $\Unom$-torsor; its isomorphism class gives us a point
\[
\ka(z) \in H^1(G(\SSS), \Un)
\,.
\]
The filtered $\phi$ variant $\ka^{F\phi}$ is defined similarly. The commutativity of the upper square is discussed at length in Dan-Cohen--Wewers \cite{CKtwo} for $n=2$; the same discussion applies here, without change. 

\segment{1711}{Discussion of $\al$}%%%%%%%%%%%%%%%%%%
We define $\al$ to be the composite
\[
\xymatrix{
X(\Zp) \ar[r]^-{\al^\m{full}} \ar[dr]_\al 	& U(X)^{F\phi} \ar[d]^{\m{projection}} \\
							& \UnFphi 
}
\]
where $\al^\m{full}$ is the unipotent Albanese map of Furusho \cite{FurushoI, FurushoII}. Our use of $u\inv$ (in place of $u$) ensures that the upper triangle commutes.

\ssegment{1112.1715}{}%%%%%%%%%%%%%%%%%%
We wish to give an explicit formula for the projection, and subsequently, an explicit formula for $\al$. Recall that $U(X)^{F\phi}$ is canonically free on two generators $x,y$ representing monodromy around $0$ and $1$ respectively 
(see for instance Furusho \cite[Notation 2.1]{FurushoII}). 
As such, it embeds as the functor of grouplike elements in the noncommutative formal power series ring space
\[
\AA^1\langle\langle x, y \rangle \rangle : R \mapsto R \langle\langle x, y \rangle \rangle.
\]
We consider the latter as a Hopf algebra for the comultiplication induced by
\begin{align*}
x \mapsto 1 \otimes x + x \otimes 1 && y \mapsto 1 \otimes y + y \otimes 1
\,.
\end{align*}
Let $\nN \subset \QQ\langle\langle x,y \rangle\rangle$ denote the Lie algebra of Lie elements, and let $V$ denote the Lie ideal with basis
\[
x, [xy], [x^2y], [x^3y], \dots
\,.
\]
We use the notational convention $[abc] = [a,[b,c]]$, etc. The claim and proposition which follow are well known; see for instance Beilinson-Deligne \cite{BeilinsonDeligne}.

\subsection*{Claim}
The representation
\[
U \to \GL(V)
\]
induced by the adjoint representation is given by
\[
\sum_w L_w w \mapsto
\begin{pmatrix}
\ddots 	&  	& 	& 	& \vdots			& \vdots \\
		& 1 			& L_x		& L_x^2/2	&L_x^3/3!		& -L_{xxxy} \\
		& 			& 1				& L_x		& L_x^2/2			& -L_{xxy} \\
		&			&				& 1				& L_x	& -L_{xy} \\
		&			&				&				& 1			& -L_y  \\
		&			&				&				&			& 1
\end{pmatrix}
\,.
\]

\begin{proof}
This is merely a verification, using the grouplike property of $\sum L_w w$ and beginning with the associated map
\[
\nN \to \fk{gl}(V)
\]
of Lie algebras. The latter is given by
\[
\sum_w l_w w \mapsto
\begin{pmatrix}
\ddots 	&  \ddots 	& 	& 	& 			& \vdots \\
		&  0			& l_x		& 	&		& -l_{xxxy} \\
		& 			& 0				& l_x		& 			& -l_{xxy} \\
		&			&				& 0				& l_x	& -l_{xy} \\
		&			&				&				& 0			& -l_y  \\
		&			&				&				&			& 0
\end{pmatrix}
\,.
\]
\end{proof}

\noindent
As a corollary, we obtain the following

\subsection*{Proposition}
\[
\al(b)=
\begin{pmatrix}
\ddots 	&  	& 	& 	& \vdots			& \vdots \\
		& 1 			& \log^p(b)		& (\log^pb)^2/2	&(\log^pb)^3/3!		& \Li_4^p(b) \\
		& 			& 1				& \log^p(b)		& (\log^pb)^2/2			& \Li_3^p(b) \\
		&			&				& 1				& \log^p(b)	& \Li_2^p(b) \\
		&			&				&				& 1			& -\log^p(1-b) \\
		&			&				&				&			& 1
\end{pmatrix}
\,.
\]

\segment{1712}{}%%%%%%%%%%%%%%
In terms of the construction summarized in the diagram of Segment \ref{1120.1}, the algorithm we hope to construct will proceed according to the following rubric.

\bigskip\noindent
\textbf{Algorithm} (sketch).
Initiate a search for integral points. This gives us a gradually increasing list $S' \subset X({\bf S})$ of points. Simultaneously with this search, we compute a $p$-adic approximation of $X(\Zp)_n$; we gradually increase $n$. This, in turn, gives us a gradually decreasing set of possible locations for integral points. The algorithm halts when (and if) all possible locations are accounted for by integral points that have already been found.

The computation of $X(\Zp)_n$ for a given $n$ and a given level of precision takes place in three steps. 

\ssegment{1712.1}{}
We search for a basis of $A(\SSSbar)_i$ for $i\le n$ among the motivic iterated integrals which are unramified over $\SSSbar$. Such a search can be made algorithmic by using, for instance, $p$-adic approximations of $p$-adic periods.

\ssegment{1712.2}{}
We compute the scheme-theoretic image of $\la$ in terms of the basis constructed in the previous step. We can then use the resulting equations to find relations between the components of $F\phi^H$, and to translate those into Coleman functions on $X(\Zp)$, without difficulty.

\ssegment{1712.3}{}
We compute $p$-adic approximations of the roots of those functions. An algorithm constructed by Besser--de Jeu for computing polylogarithms provides a method which can be generalized to arbitrary iterated integrals.

\segment{21d}{}%%%%%%
As stated in the introduction, our algorithm will have the following properties. (1) If the algorithm halts for the input $S$, then the output is equal to the set $X(\SSS)$ of solutions to the $S$-unit equation. (2) If Kim's conjecture holds for $\SSS$, and exhaustion by iterated integrals holds for an open subscheme
\[
\SSSbar \subset \SSS,
\]
then the algorithm halts for the input $S$.

\segment{21e}{Remarks}%%%%%%%%%%
Our algorithm consists of two parts which run simultaneously, neither one of which would halt individually. One part searches through the countably infinite set of possible rational points. The second part searches through a set, which is infinite (but countably so) in several ways: (1) there are infinitely many levels $n$ to consider, (2) at each level $n$, there are potentially infinitely many iterated integrals and infinitely many open subschemes $\SSSbar \subset \SSS$ to consider as we search for a basis of iterated integrals, and (3) there are infinitely many levels of $p$-adic precision which may intervene both when we attempt to construct a basis, and when we attempt to locate the zeros of the resulting Coleman functions. Nevertheless, if both conjectures hold, then after a finite amount of time, we will have a complete list of integral points, an open subscheme $\SSSbar \subset \SSS$, a basis of $A(\SSSbar)_{\le n}$ for some $n$, formulas for the Coleman functions defining $X(\Zp)_n$, and sufficient $p$-adic precision to exclude (provably) any extraneous zeros. In particular, 
we do not need a conjecture about the level $n$ at which Kim's conjecture holds in order for our algorithm to halt. We also don't need a more precise version of the exhaustion conjecture (like the one considered in segment \ref{1113.1515}), although such a conjecture would make our algorithm much more efficient, by eliminating the need to vary $\SSSbar$, and greatly limiting the set of iterated integrals that need to be considered. Finally, we expect that parts \ref{1712.1} and \ref{1712.2} of our algorithm will depend on the conjectured nonvanishing of $\ze^p(n)$ for $n$ odd $\ge 3$. However, our dependence on this conjecture is quite week: since the conjecture is known for $p$ regular, and since there are infinitely many regular primes, we can circumvent the conjecture by choosing the auxiliary prime $p$ to be regular.

%%%%%%%%%%%%%%%%%%%%%%%%%%%%%%%%%%%%%%
\section{The case $\SSS = \Spec \ZZ$}%%%%%%%%%%%%%%%%%%
%Major rev. starting 131114_1255. Backup saved to dropbox/prof13/backups

\segment{0907b}{Proposition}%%%%%%%%%%%%%%%%%%%%%%%
Assume $\ze^p(n) \neq 0$ for $n$ odd $\ge3$. Then the image of
\[
H^1\big(G(\SSS), U^n\big)_\Qp 
\xto{F\phi^H} 
H^1\big(G(\Fphi), \UnFphi\big) 
\underset{\sim}{\xto{ev_{u\inv}}} 
\UnFphi = \Qp(1) \ltimes \prod_1^n \Qp(i)
\]
consists precisely of the coordinate-plane
\[
\prod_{i \in[3,n] \mbox{ odd}} \Qp(i)
\,.
\]

\begin{proof}
For the category of mixed Tate motives over $\Spec \ZZ$, we have
\[
H^1(\QQ(1)) = H^0(\QQ(1)) = 0
\,,
\]
so the long exact sequence associated to the extension
\[
1 \to \prod_1^n \QQ(i) \to U^n \to \QQ(1) \to 1
\]
gives us an isomorphism
\[
H^1(\prod_{i=1}^n\QQ(i)) = H^1(U^n)
\,.
\]
We also have $H^1(\QQ(i)) = 0$ for $i$ even. On the other hand, for $i$ odd $\ge 3$, the map
\[
H^1(\QQ(i))\otimes \Qp \to \Qp(i)
\]
has been amply studied: the cohomology group is one-dimensional, generated by a special zeta motive $\ze^U (i)$; its image in $\Qp(i)$ is Furusho's $p$-adic special zeta value $\ze_p(i)$.
\end{proof}

\segment{}{}%%%%%%%%%%%%%%%%%%
We may now pull back along $\al$ using Proposition \ref{1112.1715} to obtain Theorem \ref{1126}. The condition $\ze^p(n) \neq 0$ for $n$ odd $\ge 3$ is known for $p$ regular and conjectured by Iwasawa theorists for $p$ arbitrary. See Examples 2.19(b) of Furusho \cite{FurushoI}.

%%%%%%%%%%%%%%%%%%%%%%%%%%%%%%%
\section{The case $S=\Spec \ZZ \setminus \{2\}$ at level $n = 4$}%%%%%%%%%%%%%%%%%
%%%%%%%%%%%%%%%%%%%%%%%%%%%%%%%%

\segment{1716}{The polylogarithmic quotient}%%%%%%%%%%%%%%%%%%%%%%
Let $U_0(X)$ denote the unipotent fundamental group of $X = \thrpl$ at the usual tangent vector at $0$. We sometimes omit either the $0$ or the $X$ from the notation.  Let $_1P_0$ denote the torsor of paths to the tangent vector $-1$ at $1$. Recall again the depth filtration $D$ of Deligne-Goncharov \cite{DelGon}, Segment 6.1: we consider the map
\[
U(X) \to \Ga
\]
induced by the natural inclusion $X \subset \Gm$, we let $D^1$ denote its kernel, and $D^n$ the descending central series of $D^1$; Deligne shows that
\[
U(X) / D^2 = \prod_1^\infty \QQ(i) \rtimes \QQ(1)
\,.
\]
We refer to $U(X) / D^2$ as the \emph{polylogarithmic quotient}.

\ssegment{}{
%The motivic iterated integrals $\ka_{i,j}$
}%%%%%%%%%%%%%%%%%%%%%%%%
As above, we let $\Mot(\SSS)$ denote the category of mixed Tate motives over $S$ with $\QQ$-coefficients, we let $U(\SSS)$ denote the unipotent radical of the fundamental group at the de Rham fiber functor, and we let $A = A(\SSS)$ denote the graded Hopf algebra of functions on $U(\SSS)$. Fix an integral point $b \in X(\SSS)$. We begin by constructing certain functions $\ka_{i,j}(b) \in A_{j-i}$ which arise naturally from the polylogarithmic quotient; we will see in a moment that these are the same as the unipotent motivic logarithms and polylogarithms defined above.
%\mar{Note that this version of motivic periods has some problems to do with Tate twists (see \cite{BrownSingle} remark 2.9)}

Let $V \subset \Lie \prod \QQ(i) \rtimes \QQ(1)$ be the Lie ideal obtained by eliminating the factor $\QQ(1)$ in the product, and let 
\[
V_b : = {_bP_0 } \overset {U_0(X)} \times V 
\]
denote the twist of $V$ by the path torsor $_bP_0$, which is a representation of the fundamental group $U_b(X)$ at $b$. The twist by $_bP_0$ does not affect the associated graded for the weight filtration, so we have
\[
\gr^W_{-2i}V_b = \QQ(i)
\,.
\]
We set 
\[
\ka_{i,j}(b) = [V_b, v_j, f_i]^{\om,\om}
\]
where $\om$ denotes the canonical fiber functor (which, we recall, is canonically equivalent to the de Rham fiber functor), and
\begin{align*}
v_j:\QQ(-j) \iso \gr^W_jV_b
&&
f_i: \gr^W_iV_b \iso \QQ(-i)
\end{align*}
are the canonical isomorphisms.

\ssegment{0913}{Coproduct formula}%%%%%%%%%%%%%%%%%%%%%%
One advantage of passing to Deligne's representation is that we obtain the following simple formula for the reduced coproduct
\[
d: A_n \to A_1 \otimes A_{n-1} + A_2 \otimes A_{n-2} + \cdots + A_{n-1}\otimes A_1
\]
for the $\ka_{i,j}$'s.

\subsection*{Lemma}
We have for $b \in X(\bf{S})$ and  $i < j \le-1$,
\[
d(\ka_{i,j}(b)) = \sum_{i<l<k} \ka_{i,l}(b) \cdot \ka_{l,j}(b)
\,.
\]

\begin{proof}
As we will see in segment \ref{0906} below, the matrix entries $\ka_{i,j}$ of the polylogarithmic quotient are the classical polylogarithms. This formula is then a special case of the Goncharov coproduct \cite[Theorem 1.2]{GonGal}; it also follows directly from the general construction of the coproduct of a \emph{framed object} in a \emph{mixed Tate category} \cite[Appendix A]{GonGal}.
\end{proof}

\ssegment{0906}{Proposition}%%%%%%%%%%%%%%%%%%%%%%%%%%%%%%
Arranging the $\ka_{i,j}$ as a matrix, we have
\[
\begin{pmatrix}
\ddots 	&				& \vdots 		\\
		& \ka_{-2,-2}(b)	& \ka_{-2,-1}(b)	\\
		&				& \ka_{-1,-1}(b)
\end{pmatrix}
=
\begin{pmatrix}
\ddots 	&  	& 	& 	& \vdots			& \vdots \\
		& 1 			& \log^U (b)		& (\log^U b)^2/2	&(\log^U b)^3/3!		& \Li^U_4 (b) \\
		& 			& 1				& \log^U (b)		& (\log^U b)^2/2			& \Li^U_3 (b) \\
		&			&				& 1				& \log^U (b)	& \Li^U_2 (b) \\
		&			&				&				& 1			& -\log^U(1-b) \\
		&			&				&				&			& 1
\end{pmatrix}
\,.
\]

\begin{proof}
 Let $c_{_b 1_0}$ denote the cocycle
\[
U(\SSS) \to U_0(X)^\dR 
\]
associated to the de Rham realization $_bP_0(X)^\dR$ of the path torsor. Let $\rho^\dR$ denote the de Rham realization of Deligne's representation
\[
U_0(X)^\dR \to \GL(V^\dR) = \GL\left( \bigoplus_{i=1}^\infty \QQ(i)^\dR \right)
\,.
\]
Then on the one hand, we have
\[
c_{_b 1_0} = \sum \Li^U_w(b)  w
\]
as grouplike elements of $A(\SSS) \langle \langle x,y \rangle \rangle $, and on the other hand,
\[
\rho^\dR \circ c_{_b 1_0} \in \GL_{A(\SSS)}(A(\SSS) \otimes V^\dR) = \GL_{A(\SSS)}\left( \bigoplus_{i=1}^\infty A(\SSS)(i) \right)
\]
is the the matrix of $\ka_{i,j}(b)$'s. So the proposition follows from the calculation of Segment \ref{1112.1715}.
\end{proof}

\ssegment{}{Variant}%%%%%%%%%%%%%%%%%%%%%%%%%%%%
The above discussion goes through with a tangent vector in place of the integral point $b$. Recall that $\ze^U(1)  = 0$, and $\ze^U(n) = 0$ for $n$ even. The result is as follows.

\subsection*{Proposition}
We have
\[
\begin{pmatrix}
\ddots 	&				& \vdots 		\\
		& \ka_{-2,-2}(1)	& \ka_{-2,-1}(1)	\\
		&				& \ka_{-1,-1}(1)
\end{pmatrix}
=
\begin{pmatrix}
\ddots 	&		&  	& 	& 	& 			& \vdots \\
		& 1		&	&	&	&			& \ze^U(5)\\
		&& 1 			& 		& 	&		& 0 \\
		&& 			& 1				& 	& 			& \ze^U(3) \\
		&&			&				& 1				& 	&  0 \\
		&&			&				&				& 1			& 0 \\
		&&			&				&				&			& 1
\end{pmatrix}
\,.
\]

%\ssegment{}{The generators $v_i$}%%%%%%%%%%%%%%%%%%%%%%
%Let
%$
%\nN(\Fphi) = \Lie U(\Fphi)
%\,.
%$
%We define elements $v_i \in \nN(\Fphi)_i$ for $i \le -1$ by 
%\[
%u\inv = \exp(\sum_{-1}^\infty v_i)
%\,.
%\]
%We then have
%\[
%\Uu U(\Fphi) = \AA^1 \langle \langle v_{-1}, v_{-2}, \dots  \rangle \rangle
%\,,
%\]
%where $\Uu$ denotes the completed universal enveloping algebra. If $w$ is a word in the $v_i$, we call the sum of the subscripts its \emph{degree} (which is equal to half its weight).

\segment{}{}%%%%%%%%%%%%%%%%%
We let $f_w$ denote the linear functional
\[
\AA^1 \langle \langle v_{-1}, v_{-2}, \dots  \rangle \rangle \to \AA^1
\]
sending a power series to the coefficient of $w$. Recall that $U(\Fphi)$ is realized as the space of grouplike elements; we continue to denote by $f_w$ the restriction of $f_w$ to the grouplike elements, which then belongs to $A(\Fphi)_n$ with $n$ equal to the degree of $w$ (we are here referring to the graded degree, equal to the sum of the subscripts, or, which is the same, to half the weight). In this notation, a basis of $A(\Fphi)_n$ is given by
\[
\set{f_w}{w \mbox{ of degree }n}
\,.
\]
To simplify notation, we write $w$ as a word in $\NN$, with ``$.$'' denoting multiplication in the power-series ring; so for instance,
\[
f_{1^{.2}} = f_{1.1} = f_{v_{-1}v_{-1}}
\,.
\]

%\ssegment{}{The $p$-adic period}%%%%%%%%%%%%%%%%%%%%%%
%We denote $F\phi$-realization
%\[
%A \to A(\Fphi)
%\]
%by $F\phi$, placing the symbol $F\phi$ sometimes on the left, as in $F\phi(\Li^U_w b)$, and sometimes as an exponent, as in $\Li_w^{F\phi}(b)$ or $(\Li^U_w b)^{F\phi}$ (there is no ambiguity by Segment \ref{4_14_f}). We have
%\[
%(\Li^U_w b)^{F\phi}(u\inv) = \Li_w^p b
%\,,
%\]
%the $p$-adic polylogarithmic value of Furusho \cite{FurushoI}. 

\segment{1119}{Determination of bases; computation of $F\phi$}%%%%%%%%%%%%%%%%%%%%%%%%%%%
We now begin to work our way up the latter, developing more tools as they become necessary. In Segment \ref{1119}.n (n=1,2,3,4) we find a basis of $A(\SSS)_n$ ($\SSS = \Spec \ZZ \setminus \{2\}$). For $n = 1,2,3$ we also compute the matrix associated to
\[
F\phi: A(\SSS)_n \to A(\Fphi)_n
\,.
\]
Let us abbreviate
\[
A(\SSS) =: A.
\]

\ssegment{}{}%%%%%%%%%%%%%%%%%%%%%%%%%%%%%%
A basis of $A_1$ is given by $\{\log^U 2 \} $. A basis of $A(\Fphi)_1$ is given by $\{f_1\}$. The $1 \times 1$ matrix associated to $F\phi_1$ is given by
\[
(\log^U 2)^{F\phi}(v_{-1}) = \log^p(2)
\,.
\]

\ssegment{}{}%%%%%%%%%%%%%%ppp

%\sssegment{}{}%%%%%%%%%%%%%%%%%%%%%
A basis of $A_2$ is given by $(\log^U 2)^2$. A basis of $A(\Fphi)_2$ is given by $\{f_{1^{.2}}, f_2\}$. Using the shuffle product of Segment \ref{0907}, we find that
\[
(\log^{F\phi} 2)^2(v_{-1}^2)
= 2((\log^{F\phi} 2)(v_{-1}))^2 = 2(\log^p 2)^2
\,,
\]
while
\[
(\log^{F\phi} 2)^2(v_{-2})=0
\,,
\]
from which we conclude that the $1 \times 2$ matrix associated to $(F\phi)_2$ is given by
\[
\begin{pmatrix}
2(\log^p 2)^2
\\
0
\end{pmatrix}
\,.
\]

\sssegment{1009}{Identity for $\Li^U_2(1/2)$}%%%%%%%%%%%%%%%%%
If $b \in X(\ZSinv)$, then $\Li^U_2 b = r(\log^U 2)^2$ with $r \in \QQ$. We compute $r$ in two ways. Our first method is to evaluate (the $F\phi$ realization) at $v_{-1}^2$. Since
\[
v_{-1}(V_{3,b})=
\begin{pmatrix}
0 & \log^p b & 0 \\
& 0 & -\log^p(1-b) \\
&& 0
\end{pmatrix}
\]
we find that $(\Li^U_2b)(v_{-1}^2) = (v_{-1}(V_{3,b}))^2_{1,3} = -(\log^pb)\log^p(1-b)$, so
\[
r = -\frac{(\log^p b)\log^p(1-b)}{2(\log^p 2)^2}
\,.
\]
Our second method is to apply the boundary map
\[
A_2 \to A_1\otimes A_1
\,.
\]
We apply the comultiplication map to the $(3,2,1)$-framed object $V_{3,b}$ with associated $k$-framed objects $\ka_{i,j}$ with $k=j-i = 0,1,2$ given by
\[
\begin{pmatrix}
1 & \log^U b & \Li^U_2 b \\
& 1 & - \log^U(1-b) \\
&& 1
\end{pmatrix}
\]
and eliminate the edge-terms in $A_0 \otimes A_2 + A_2\otimes A_0$ to obtain
\[
d(\Li^U_2b)= - (\log^U b)\otimes \log^U(1-b)
\,.
\]
On the other hand, we have
\[
d((\log^U 2)^2) = 2 (\log^U 2)\otimes (\log^U 2)
\,,
\]
so we again obtain
\[
-(\log^U b) \otimes (\log^U(1-b)) = r \cdot 2(\log^U 2) \otimes (\log^U 2)
\,,
\]
as expected. For instance, for $b=1/2$, we get a well-known identity.

\subsection*{Lemma}
\[
\Li^U_2(1/2) = -\frac{1}{2}(\log^U 2)^2
\,.
\]

\ssegment{}{}%%%%%%%%%%%%%%%%%%%%%%%%%%%
We begin our discussion of the level $n=3$.

\sssegment{0907c}{}%%%%%%%%%%%%%%
Let us record in detail how we go about computing the value $\Li^{F\phi}_3(b)(w)$, for
\[
w \in \{ v_{-1}^3, v_{-1}v_{-2}, v_{-2}v_{-1}, v_{-3} \}
\]
a word of degree $3$. The lower-right $4 \times 4$ block
\[
M = 
\begin{pmatrix}
	 1				& \log^U b		& (\log^U b)^2/2			& \Li^U_3 b \\
					& 1				& \log^U b	& \Li^U_2 b \\
				&				& 1			& -\log^U(1-b) \\
					&				&			& 1
\end{pmatrix}
\] 
of the matrix appearing in Proposition \ref{0906} reminds us that
\[
\Li^U_3(b) = \ka_{-4,-1}(b)
\]
may be represented by the $(-4,-1)^\m{st}$ matrix entry of the representation $(V_4)_b$, where $V_4$ denotes the quotient
\[
\bigoplus_1^4\QQ(i)
\]
of $V$, and the basis element of $\QQ(i)_\dR$ is taken to have degree $-i$ as usual. The action $v_{-i}((V_4)_b)$ of $v_{-i}$ is given by the $i^{\rm{th}}$ superdiagonal of $\log (M^{F\phi}(u\inv))$. We find that
\begin{align*}
v_{-1}(V_{4,b}) &=
\begin{pmatrix}
0& \log^pb && \\
& 0 & \log^pb &\\
&& 0 & -\log^p(1-b) \\
&&& 0
\end{pmatrix}
\\
v_{-2}(V_{4,b}) &=
\begin{pmatrix}
0 &0 & 0 & \\
&0&0& \frac{1}{2}(\log^pb)\log^p(1-b) + \Li^p_2 b \\
&&0&0\\
&&&0
\end{pmatrix}
\\
v_{-3}(V_{4,b}) &=
\begin{pmatrix}
0 &0 & 0 &  -\frac{1}{12} (\log^p b)^2\log^p(1-b) - \frac{1}{2}(\log^p b)\Li_2^p b + \Li_3^p b \\
&0&0& 0 \\
&&0&0\\
&&&0
\end{pmatrix}
\end{align*}
$\Li_3^{F\phi}(b)(w)$ is given by the $(-4,-1)^{\rm {st}}$ matrix entry of $w(V_{4,b})$. So multiplying the above matrices and reading off the north-eastern corner, we find
\begin{align*}
\Li_3^{F\phi}(b)(v_{-1}^3) &= -(\log^p b)^2 \log^p(1-b)
\\
\Li_3^{F\phi}(b)(v_{-1}v_{-2}) &= \frac{1}{2}(\log^p b)^2\log^p(1-b) +(\log^p b)\Li_2^p b
\\
\Li_3^{F\phi}(b)(v_{-2}v_{-1}) &= 0
\\
\Li_3^{F\phi}(b)(v_{-3}) &= -\frac{1}{12} (\log^p b)^2\log^p(1-b) - \frac{1}{2}(\log^p b)\Li_2^p b + \Li_3^p b 
\,.
\end{align*}

\sssegment{}{Proposition}%%%%%%%%%%%%%%%%%%%%%%%%%
A basis of $A_3$ is given by $ \{ (\log^U 2)^3, \ze^U(3) \} $.

\begin{proof}
We use the complex $A(3)$ (\ref{1120.2}) to conclude that $\dim A_3 =2$. Recall that for $i \ge 2$, the Ext groups $\Ext^1(\QQ(0), \QQ(i))$ for mixed Tate motives unramified over an open subscheme $\SSS \subset \Spec \ZZ$ don't depend on $\SSS$, and that in particular, for $i$ odd $\ge 2$, the Ext group is generated by $\ze^U(i)$. On the other hand, denoting the coproduct by $\nu$, we have
\begin{align*}
\nu\big( (\log^U 2)^3 \big) &= \big( \nu(\log^U 2) \big)^3
\\
	&= \big( 1 \otimes (\log^U 2) + (\log^U 2) \otimes 1 \big)^3
\\
	&= 1 \otimes (\log^U 2)^3 + 3(\log^U 2)\otimes (\log^U 2)^2 + 3(\log^U 2)^2\otimes (\log^U 2) + (\log^U 2)^3 \otimes 1
\,,
\end{align*}
so
\begin{align*}
d\big( (\log^U 2)^3 \big) &= 3(\log^U 2)\otimes (\log^U 2)^2 + 3(\log^U 2)^2\otimes (\log^U 2)
\\
	& \neq 0
\,.
\end{align*}
Since
\[
\ker d = H^1(A(3)) = \Ext^1(\QQ(0), \QQ(3))
\,,
\]
the proposition follows.
\end{proof}

\sssegment{}{}%%%%%%%%%%%%%%%%%%%
A basis of $A(\Fphi)_3$ is given by $ \{ f_{1^{.3}}, f_{1.2}, f_{2.1}, f_{3} \} $. Applying the above computations, we find that the matrix of $(F\phi)_3$ is given by
\[
\begin{pmatrix}
(\log^p 2)^3 & 0 \\
0&0 \\
0&0\\
0& \ze^p(3)
\end{pmatrix}
\,.
\]

\sssegment{1737}{Identity for $\Li^U_3(1/2)$}%%%%%%%%%%%%%%%%%%%%%%%%%

Let $b \in X(S)$ be an $S$-integral point ($S=\Spec \ZZ \setminus \{2\}$). The computations of segments \ref{1009}, \ref{0907c} are also useful for expanding $\Li^U_3 (b)$ as a linear combination
\[
\Li^U_3(b) = s\ze^U(3) + t(\log^U 2)^3
\]
of our basis elements. We should emphasize at this point, that our goal here is not to derive exact expressions for $s$ and $t$, but rather to develop methods for computing $p$-adic approximations which may be applied algorithmically. It is also important that our approximations be $p$-adic rather than complex, since our algorithm proceeds by computing increasingly good $p$-adic approximations of certain functions whose coefficients depend on the computation being performed here. For this reason we have to use $p$-adic periods and not complex periods.

Evaluating at $v_{-1}^3$, we find that
\[
t= -\frac{(\log^p b)^2 \log^p(1-b)}{6(\log^p 2)^3}
\,.
\]
Evaluating at $v_{-3}$, we find that
\begin{align*}
s 	&= \frac{ - \frac{1}{12} (\log^p b)^2 \log^p(1-b) - \frac{1}{2} (\log^p b) \Li_2^p b + \Li_3^p b }{\ze^p(3)}
\\
	&= \frac{\frac{1}{6} (\log^p b)^2\log^p(1-b) + \Li_3^p b }{ \ze^p(3) } 
	\,.
\end{align*}
For instance, for $b=1/2$, we have
\[
\Li^U_3(1/2) = \frac{-\frac{1}{6} (\log^p 2)^3 + \Li^p_3(1/2)  }{ \ze^p(3)} \ze^U(3) + \frac{1}{6} (\log^U 2)^3
\,.
\]

%\sssegment{}{}%%%%%%%%%%%%%%%%%%%%%%%
In particular, we've shown that
\[
\frac{-\frac{1}{6} (\log^p 2)^3 + \Li^p_3(1/2)  }{ \ze^p(3)} 
\,,
\]
a priori $\in \Qp$, is actually $\in \QQ$. Computations based on \textit{Lip service} \cite{Lip} show that it is $p$-adically close to $7/8$ for several small primes $p$. We will denote this number by $\widetilde{7/8}$. In this notation, we have
\[
\Li^U_3(1/2) = \widetilde { (7/8) }\ze^U(3) + \frac{1}{6}(\log^U 2)^3
\]
on the motivic level. Two proofs that this identity holds with
\[
 \widetilde { (7/8) } = 7/8,
\]
one given by Hidekazu Furusho using $p$-adic polylogarithms, and one given to us by one of the referees using complex polylogarithms, appear in the appendices.

\ssegment{1112}{Proposition}%%%%%%%%%%%%%%%%%%
The set $\Bb=\{(\log^U 2)^4, (\log^U 2)\ze^U(3), \Li^U_4(1/2)\}$ forms a basis of $A_4$.

\begin{proof}
A glance at the complex $A(4)$ reveals that the dimension is $3$. By the computations above, a basis for
\[
A_1 \otimes A_3 + A_2 \otimes A_2 + A_3 \otimes A_1
\]
is given by
\[
\{ (\log^U 2)\otimes \ze^U(3), (\log^U 2)\otimes (\log^U 2)^3, (\log^U 2)^2\otimes (\log^U 2)^2, (\log^U 2)^3\otimes (\log^U 2), \ze^U(3)\otimes(\log^U 2) \}
\,.
\]
Using Lemma \ref{0913} together with Proposition \ref{0906} to compute reduced coproducts
\[
d: A_4 \to A_1 \otimes A_3 + A_2 \otimes A_2 + A_3 \otimes A_1
\]
we find that 
\[
d(\Bb) = 
\begin{pmatrix}
0 & 1& -\widetilde{7/8} \\
4 &0 & -1/6 \\
6&0&-1/4 \\
4&0&-1/6 \\
0&1&0
\end{pmatrix} \,.
\tag{M}
\]
We can verify numerically that
\[
\widetilde{7/8} \neq 0
\,,
\]
so the proposition follows.
\end{proof}

\ssegment{1112.12}{Summary}%%%%%%%%%%%%%%%%%%%%%%%%%%
We've found that a basis of
\[
A(\SSS)_1 \times \prod_1^4 A(\SSS)_i
\]
is given by
\[
\{(\log^U 2)', \log^U 2, (\log^U 2)^2, (\log^U 2)^3, \ze^U(3), (\log^U 2)^4, (\log^U 2)\ze^U(3), \Li^U_4(1/2)\} \tag{C}
\,,
\]
where $(\log^U 2)'$ denotes $\log^U 2$ regarded as an element of the left-hand copy of $A_1$. 

\segment{}{Computation of image of $\la$}%%%%%%%%
Our next task is to compute the scheme-theoretic image of the map
\[
\la: H^1(G(\SSS), \Un) \to A(\SSS)_1 \times \prod_1^4 A(\SSS)_i
\]
in terms of these coordinates. 

\ssegment{1112.1207}{}%%%%%%%%%%%%%%%%%%%%%%%%
We begin by considering arbitrary generators $\nu_i$, $i$ odd $\le -1$, of $U(\SSS)$. These arbitrary generators give us a second, \textit{abstract} basis of the target space. We write $\phi_{1^{.2}}$ for the element of $A(\SSS)$ dual to $(\nu_{-1})^2$, and similarly for any word in the $\nu_i$; in this notation, the abstract basis is given by
\[
\{ \phi_1', \phi_1, \phi_{1^{.2}}, \phi_{1^{.3}}, \phi_3, \phi_{1^{.4}}, \phi_{1.3}, \phi_{3.1} \}
\tag{A}
\,.
\]
We wish to compute the image of $\la$ in terms of this basis.

\ssegment{r_15_3}{}%%%%%%%%%%%%%%%%%%
We use the notation
\[
\la = \la_1' \times \la_1 \times \la_2 \times \la_3 \times \la_4
\]
for the components of $\la$. We write $r^D = \Lie \rho^D$. An arbitrary element $\rho$ of
\[
\Hom^\Gm(U(\SSS), U^4)
\]
with corresponding element $r$ of
\[
\Hom^{\gr}(\nN(\SSS), \nN^4)
\]
may be specified by setting
\begin{align*}
r^D r(v_{-1}) = 
\begin{pmatrix}
0&a\\
&0&a \\
&&0&a\\
&&&0&-b\\
&&&&0
\end{pmatrix},
&&
r^Dr(v_{-3}) = 
\begin{pmatrix}
0&&&&0\\
&0&&&-d\\
&&0&&0\\
&&&0&0\\
&&&&0
\end{pmatrix},
\end{align*}
with $a,b,d, \in \QQ$. We wish to solve
\[
\la_1'(\rho) = x_1' \phi_1'
\]
for $x_1' \in \QQ$. In the following diagram, $\la_1'(\rho)$ appears as the composite from upper left to the far right.
\[
\xymatrix{
U(\SSS) \ar[r]^-\rho \suphook[d] & 
U^4 \ar[r]^-{\rho^D} \suphook[d] & 
U^{5\times 5} \suphook[d]
\\
\Uu\nN(\SSS) \ar[r]^-R & 
\Uu\nN^4 \ar[r]^-{R^D} & 
\End  \bigoplus_1^5 \QQ(i) \ar[r]_-{\pi_{-2,-3}} & 
\AA^1
\\
\nN(\SSS) \ar[r]_-r \suphook[u] & 
\nN^4 \ar[r]_-{r^D} \suphook[u] &
\nN^{5\times5} \suphook[u]
}
\]
We need to \textit{evaluate} it on the element $v_{-1} \in \nN(\SSS)$ on the lower left, not a priori in its domain of definition. To this end we consider the induced maps of  enveloping algebras, as in the central row. We find that the central composite $\pi_{-2,-3} \circ R^D \circ R$ is the unique linear functional on $\Uu\nN(\SSS)$ restricting to $\la_1'(\rho)$ on $U(\SSS)$, and so we have
\[
x_1' = (\pi_{-2,-3} \circ R^D \circ R)(v_{-1}) = (\pi_{-2,-3} \circ r^D \circ r)(v_{-1}) = a.
\]
Similarly, by computing appropriate words in the two matrices above and extracting appropriate matrix entries, we find that
\begin{align*}
\la(\rho) = a \phi_1' -b\phi_1 -ab \phi_{1^{.2}} 	&- a^2b\phi_{1^{.3}} - d \phi_3 -a^3b\phi_{1^{.4}} - ad \phi_{1.3} 
\end{align*}
If we denote an arbitrary element of $A_1 \times \prod_1^4 A_i$ by
\[
x_1'\phi_1' + x_1\phi_1+x_{1^{.2}}\phi_{1^{.2}}+ \cdots + x_{3.1}\phi_{3.1}
\]
then we've found that equations for the image of $\la$ are as follows.
\begin{align}
x_{1^{.2}} &= x_1x_1' \\
x_{1^{.3}} &= x_1{x_1'}^2 \\
x_{1^{.4}} &= x_1{x_1'}^3 \\
x_{1.3} &= x_1'x_3 \\
x_{3.1} &=0
\end{align}

\ssegment{1736}{}%%%%%%%%%%%%%%%%
Our next task is to rewrite equations \ref{1112.1207}.(1-5) in terms of the concrete basis \ref{1112.12}.(C), allowing ourselves to impose certain conditions on our arbitrary generators $\nu_i$ of $U(\SSS).
$\footnote{We have learned that the computation which follows is quite similar to the algorithm constructed by Francis Brown in \textit{On the decomposition of motivic multiple zeta values} \cite{BrownDecomp}.} 
Since $\dim A_1 = 1$, $\log^U 2$ is a scalar multiple of $\phi_1$, so after possibly replacing $\nu_{-1}$ by a scalar multiple, we may assume
\begin{align*}
(\log^U 2)'	&= \phi'_1 \tag1 \\
\log^U 2 	&= \phi_1 \tag2 \,.\\
\end{align*}
The formula $(\phi_1)^n = n!\phi_{1^{.n}}$ then implies
\begin{align*}
(\log^U 2)^2 &= 2\phi_{1^{.2}} \tag3 \\
(\log^U 2)^3 &= 6\phi_{1^{.3}} \,. \tag4 \\
\end{align*}
Similarly, $\phi_3$ and $\ze^U(3)$ belong the the one-dimensional subspace $\Ext^1(\QQ(0), \QQ(3))$ of $A_3$, so after possibly replacing $\nu_{-3}$ by a scalar multiple, we may assume
\[
\ze^U(3) = \phi_3 \tag{5} \,.
\]
Subsequently, we have
\begin{align*}
(\log^U 2)^4 &= 24\phi_{1^{.4}} \tag6 \\
(\log^U 2)\ze^U(3) &= \phi_{1.3}+\phi_{3.1} \tag7 \,.
\end{align*}
We use our computation \ref{1112}.(M) of the injection
\[
d : A_4 \to A_1\otimes A_3 + A_2\otimes A_2 + A_3\otimes A_1
\]
to expand $\Li^U_4(1/2)$ in terms of the abstract basis:
\begin{align*}
d(\Li^U_4(1/2)) & = -\widetilde{7/8}(\log^U 2)\otimes \ze^U(3) - \frac{1}{6} (\log^U 2)\otimes(\log^U 2)^3-\frac{1}{4}(\log^U 2)^2\otimes(\log^U 2)^2 - \frac{1}{6}(\log^U 2)^3\otimes \log^U 2
\\
	&= -\widetilde{7/8}\phi_1\otimes \phi_3 - \phi_1\otimes\phi_{1^{.3}} - \phi_{1^{.2}}\otimes\phi_{1^{.2}} - \phi_{1^{.3}} \otimes \phi_1
\\
	&= -\widetilde{7/8}d(\phi_{1.3}) - d(\phi_{1^{.4}})
	\,.
\end{align*}
so
\[
\Li^U_4(1/2) = -\widetilde{7/8}\phi_{1.3} - \phi_{1^{.4}}
\tag{8}
\,.
\]

\ssegment{1112.1654}{}%%%%%%%%%%%%%%%%%%%%%%
If we denote an arbitrary linear combination of our concrete basis elements \ref{1112.12}.(C) by
\begin{align*}
y=
y'_l (\log^U 2)'+ y_l\log^U 2 + y_{l^2} &(\log^U 2)^2 + y_{l^3}(\log^U 2)^3 + y_\ze \ze^U(3) 
\\
	& + y_{l^4}(\log^U 2)^4 + y_{l\ze}(\log^U 2)\ze^U(3) + y_L Li_4(1/2)
\,,
\end{align*}
then in terms of these coordinates, equations \ref{1112.1207}.(1-5) become
\begin{align*}
2 y_{l^2} &= y_ly_l' \tag1 
\\
6 y_{l^3} &= y_l{y_l'}^2 \tag2
\\
24 y_{l^4} - y_L &= y_l{y_l'}^3 \tag3
\\
- \widetilde{7/8}y_L &= y_l' y_\ze \tag4
\\
y_{l\ze} &=0 \tag5
\,.
\end{align*}

\segment{1738}{Poof of Theorem \ref{1120.14}}%%%%%%%%%%%%%%%%%%%%%%%
The map
\[
A(\SSS)_1 \times \prod_1^4 A(\SSS)_i \to A(\Fphi)_1 \times \prod_1^4 A(\Fphi)_i \xto{ev_{u\inv}} {\AA^5}
\]
sends $y$ to
\[
\begin{pmatrix}
(\log^p2)y_l' \\
(\log^p2)y_l \\
(\log^p2)^2 y_{l^2} \\
(\log^p2)^3 y_{l^3} + \ze^p(3) y_\ze \\
(\log^p2)^4 y_{l^4} + (\log^p2)\ze^U(3)y_{l\ze} + \Li^p_4(1/2) y_L
\end{pmatrix}
\]
which, according to equations \ref{1112.1654}.(1-5), equals
\[
\begin{pmatrix}
(\log^p2)y_l' \\
(\log^p2)y_l \\
\frac{(\log^p2)^2}{2} y_l y_l' \\
\frac{(\log^p2)^3}{6} y_l{y_l'}^2 + \ze^p(3)y_\ze \\
\frac{(\log^p2)^4}{24} y_l{y_l'}^3 
- \widetilde{8/7} \left( \frac{(\log^p2)^4}{24} + \Li^p_4(1/2) \right) y_l' y_\ze
\end{pmatrix}
\,.
\]
Denoting the standard coordinates of $\AA^5$ by $X_1, Y_1, Y_2, Y_3, Y_4$, we obtain the two equations
\begin{align*}
Y_2 &= \frac{1}{2} X_1Y_1 \\
Y_4 &= 
\frac{X_1^3Y_1}{24} - \widetilde{8/7} 
\left( 
		\frac{(\log^p2)^3}{24\ze^p(3)} + \frac{\Li^p_4(1/2)}{(\log^p2)\ze^p(3)}
\right)
\left(
X_1Y_3 - \frac{X_1^3Y_1}{6}
\right) 
\,.
\end{align*}
We pull back along
\[
X(\Zp) \xto{\al} U^{4,F\phi} \inj U^{5\times 5, F\phi} \xto{pr} \AA^5 
\]
using Proposition \ref{1112.1715}, to obtain the two functions given in the theorem.

%%%%%%%%%%%%%%%%%%%%%%%%%%%%%%%%%%%%
\section{Experimental verification of Kim's conjecture for $S=\Spec \ZZ \setminus \{2\}$}%%%
\label{r_15_4}%%%%%%%%%%%%%%%%%%%%%%%%%%%%%
%%%%%%%%%%%%%%%%%%%%%%%%%%%%%%%%%%%%%%

\segment{}{}%%%%%%%%%%%%%
It is an easy exercise to show that the set $X(S)$ of $S$-integral points of
$X=\PP^1\setminus\{0,1,\infty\}$ consists of the three points $1/2,2,-1$. Kim's
conjecture predicts that for $n$ sufficiently large we have
\[
   X(S) = \alpha^{-1}\Big(F\phi\big(H^1(G(\bS),U^n)\big)\Big)\subset X(\ZZ_p).
\]
We are able to verify this experimentally for $n=4$ and for small
primes $p$. 

\segment{}{}%%%%%%%%%%%%%%%%
We use the computer algebra system {\tt sage}. The code we have written
consists of two files, {\em localanalytic.sage} and {\em Lip.sage} \cite{LocLip}. The file {\em
  localanalytic.sage} contains the definition of a sage class {\em
  pAdicLocalAnalyticFunction}. An object in this class represents a locally
analytic function on an open subset of $\AA^1(\ZZ_p)=\ZZ_p$. In plainer words,
this means that the function is defined on a subset of  all residue disks, and
on each such residue disk it is given by a formal power series over $\QQ_p$.
The main feature of these objects is a routine that determines the (finite)
set of zeroes of a given function. The file {\em Lip.sage} contains an
implementation of the algorithm described in \cite{Lip} which creates
polylogarithmic functions ${\rm Li}_k^p$, for $k=1,..,n$.

\segment{}{}%%%%%%%%%%%%%%%%%%%
Here is a sample session where we verify Kim's conjecture for $S=\{2\}$, $n=4$
and $p=11$. The commands
\begin{sageverbatim}
  runfile localanalytic.sage
  runfile Lip.sage
\end{sageverbatim}
load the code. With
\begin{sageblock}
  p=11
  LIP=make_Li(p,'all',4,30)
\end{sageblock}
we construct the functions ${\rm Li}_k^p$ for $k=1,\ldots,4$, which are stored
in the list LIP. The last parameter is the relative precision used for
computations in $\QQ_p$. As a side effect, we also have commands available to
compute $\log^p(x)$ for any $p$-adic integer $x$ and $\zeta^p(k)$ for
$k=2,\ldots,4$.  For example:
\begin{sagecommandline}
  sage: logp(2^3)-3*logp(2)
  sage: zetap(3)
  sage: zetap(4)
\end{sagecommandline}
Now we check that the rational constant $\widetilde{(7/8)}$ from Segment
5.2.3.4 is $p$-adically close to $7/8$:
\begin{sagecommandline}
  sage: c1=(-logp(2)^3/6+Lip(3,1/2))/zetap(3)
  sage: c1-7/8
\end{sagecommandline}
We create the two functions $F_1$ and $F_2$ from Segment 5.4.
\begin{sagecommandline}
  sage: R.<x1,y1,y2,y3,y4>=Qp(p,30)[]
  sage: c2=c1^(-1)*(logp(2)^3/(24*zetap(3))+Lip(4,1/2)/(logp(2)*zetap(3)))
  sage: F1=subs_pol(y2-1/2*x1*y1,LIP)
  sage: F2=subs_pol(y4-x1^3*y1/24+c2*(x1*y3-x1^3*y1/6),LIP)
\end{sagecommandline}
We compute the set of zeroes of both functions:
\begin{sagecommandline}
  sage: Z1=F1.zeroes()
  sage: len(Z1)
  sage: Z2=F2.zeroes()
  sage: len(Z2)
\end{sagecommandline}
Finally, we check that the set of common zeroes is exactly the set of
$S$-integral points:
\begin{sagecommandline}
  sage: Z=padic_common_zeroes([F1,F2],20)
  sage: [z.rational_reconstruction() for z in Z]
\end{sagecommandline}

\segment{}{}%%%%%%%%%%%%%%%%%
We have run the previous commands for all primes $p=3,5,\ldots,29$ and have always
got the same results:
\begin{itemize}
\item[(a)]
  The constant $\widetilde{(7/8)}$ is $p$-adically as close to $7/8$ as can be
  expected from the chosen precision.
\item[(b)]
  The two functions $F_1$, $F_2$ have exactly the set $\{2,1/2,-1\}$ of common
  zeroes. 
\end{itemize}
Even though (a) gives very convincing evidence that
$\widetilde{(7/8)}=7/8$, it does not and cannot provide a proof. In contrast,
(b) can in principle be rigourously proved with our methods. In order to do
this, we would have to make sure that the results of our calculation hold up
to the stated precision (for instance, that $c_1\equiv 7/8\pmod{11^{24}}$ in
line 8-9). The point is that we can determine the exact number of zeroes of a
$p$-adic analytic function from a sufficiently close approximation. Since we
know by Theorem \ref{} that the $S$-integral points are zeroes of $F_1$ and
$F_2$, we can hope to prove that there are no more. 

%%%%%%%%%%%%%%%%%%%%%%%%%%%%%%%%
\appendix%%%%%%%%%%%%%%%%
%%%%%%%%%%%%%%%%%%%%%

%%%%%%%%%%%%%%%%%%%%%%%%%%%%%%%%%%%%%%%%%
\section{$\widetilde{7/8} = 7/8$ via $p$-adic polylogarithms}%%
%%%%%%%%%%%%%%%%%%%%%%%%%%%%%%%%%%%%%%%%

\segment{r28a}{}
In this section we give a proof of the relation
\[
      \Li^p_3(1/2)=\frac{7}{8}\zeta^p(3)+\frac{1}{6}\log^p(2)^3
      \,.
\]
As an immediate corollary, we obtain:
\[
\widetilde{7/8} = 7/8 .
\]
In other words, the motivic identity follows from the corresponding identity of $p$-adic periods. See segment \ref{1737}. The proof, which was suggested to us by H. Furusho, is based on the proof of its complex cousin:
\[
     \Li^\infty_3(1/2)=\frac{7}{8}\zeta^\infty(3)+\frac{1}{6}\log^\infty(2)^3-
         \frac{1}{2}\log^\infty(2)\zeta^\infty(2),
\tag{$*$}
\]
which is proved, for instance, in segment 6.12 of Lewin \cite{Lewin}. 

Proving identities as above over the $p$-adics is actually easier than over
the complex numbers, because the polylogarithmic functions are, unlike in the
complex case, single valued. There is, however, a subtle point having to do
with the value of $\Li^p_n(z)$ at $z=1$ (which is, by definition,
$\zeta^p(n)$). Explaining this argument in detail is essentially the main point of
this section. Our main reference is \cite{FurushoI}. For more on the context of identities of $p$-adic polylogarithms, the reader may consult, for instance, Coleman \cite{ColemanDilogs} or Wojtkowiak \cite{WojtkowiakNote}.

\segment{}{}
We remark that Furusho's proof does not depend on the conjectural nonvanishing of $\ze^p(3)$. Indeed, we may start by choosing a prime $p$ for which the nonvanishing is known (e.g. any regular prime). Alternatively, for any given prime $p$ we may simply verify the nonvanishing computationally. Since the identity of $p$-adic polylogarithms implies the corresponding identity of motivic polylogarithms and vice versa, the result for arbitrary $p$ then follows from the result for the given prime $p$.

\segment{}{}
Francis Brown has pointed out to us that the same result could also be deduced directly from the complex identity via the construction of \cite{BrownSingle}. A second proof, using this method, and spelled out in detail by one of the referees, is given in Appendix B.

\segment{Lidiffeq}{}%%%%%%%%%%%%%%%%%%%%
Let $\CC_p$ denote the completion of an algebraic closure of $\QQ_p$. For the remainder of this appendix, we shed the superscript $p$ above our $p$-adic polylogarithms. On the
open unit disk, the polylogarithmic functions can be easily defined as 
convergent power series. In fact, one sets
\[
%\label{eq5}
      \Li_n(z) := \sum_{k\geq 1} \frac{z^k}{k^n},
\]
for $n\geq 1$ and $z\in\CC_p$, $\abs{z}<1$. It is easy to see that the series
converges uniformly on every closed disk $\abs{z}\leq r$ with radius
$r<1$. Therefore, the functions $\Li_n(z)$ are analytic on the open unit
disk. They satisfy the differential equation
\[
%\label{Lidiffeq}
  \frac{d}{dz} \Li_n(z) = \begin{cases} \displaystyle
        \quad\frac{1}{z}\,\Li_{n-1}(z), & n\geq 2, \\ \displaystyle
        \quad \frac{1}{1-z},           & n=1.
                          \end{cases}
\tag{$*$}
\]

\segment{}{}%%%%%%%%%%%%%%%
Using Coleman's theory of $p$-adic integration, one can extend $\Li_n$ to
functions on $Y:=\PP^1(\CC_p)-\{0,1,\infty\}$, see \cite{FurushoI}. To do this we
choose an arbitrary element $a\in\CC_p$ and let $\log^a:\CC_p^\times\to\CC_p$
denote the branch of the $p$-adic logarithm defined by $\log^a(p):=a$. Then
there exist unique functions
\[
    \Li_n^a:Y\to\CC_p
\]
with the following properties.
\begin{itemize}
\item[(a)]
  $\Li_n^a$ is a {\em Coleman function} in the sense of \cite{Besser} and
  \cite{FurushoI}. (Note that the definition of Coleman functions depends on the
  choice of the branch $\log^a$ of the $p$-adic logarithm.)
\item[(b)]
  $\Li_n^a$ satisfies the differential equation
  \[
    \frac{d}{dz} \Li_n^a(z) = \begin{cases}
          \frac{1}{z}\Li_{n-1}^a(z), & n\geq 2, \\
          \frac{1}{1-z},             & n=1.
                             \end{cases}
  \]
\item[(c)]
  On the open unit disk $\abs{z}<1$, $\Li^a(z)$ agrees with $\Li_n(z)=\sum_k
  z^k/k^n$. 
\end{itemize}
For instance, 
\[
      \Li_1^a(z)=-\log^a(1-z).
\]

\segment{relrem1}{Remark}%%%%%%%%%%%
 Let 
  \[
       V:=\set{z\in Y} {\abs{z},\abs{1-z},\abs{z^{-1}}=1}\subset Y
  \]
  be the complement of the open unit disks around $0,1,\infty$. The
  restriction of the functions $\log^a(z)$ and $\Li_n^a(z)$ to $V$ are
  independent of the chosen branch of the logarithm. For a point $z\in
  V$ we may therefore write $\log(z)$ and $\Li_n(z)$. 
  
  \segment{}{}%%%%%%%%%%%%%%%%%%%%%
  In order to define $p$-adic zeta values, Furusho has shown that the functions
$\Li_n^a(z)$ have a well defined value for $z=1$ if $n\geq 2$. We need to
recall what this means. Let $x\in\{0,1,\infty\}$ and $z_x$ be a local
parameter at $x$ (e.g.\ $z_0:=z$, $z_1:=1-z$, $z_\infty=z^{-1}$). If $f$ is a
Coleman function on $Y$ then there exists a radius $0<r<1$
such that $f$ is defined on the open admissible subset 
\[
      U_x(r):=\set {z_x\in\CC_p} {r\leq\abs{z_x}<1} \subset Y.
\]
Furthermore, $f$ admits an expansion of the form
\[
    f(z_x) = f_0(z_x)+f_1(z_x)\log^a(z_x)+\ldots+f_m(z_x)\log^a(z_x)^m,
\]
where $f_0,\ldots,f_m$ are analytic functions on $U_x(r)$. We say that $f$ is
{\em defined} at $z=x$ if the functions $f_0,\ldots,f_m$ extend to analytic
functions on the whole open disk
\[
        U_x =\set{ z_x\in\CC_p } {\abs{z_x}<1 }
\]
and, moreover, $f_i(x)=0$ for $i=1,\ldots,m$. If this is the case, then $f(x):=
f_0(x)$ is called the {\em value} of $f$ at $z=x$. 

It is clear that the sum $f=g+h$ of two Coleman functions is defined at $x$ if
both $g$ and $h$ are defined at $x$. Moreover, if this is the case then
$f(x)=g(x)+h(x)$. Furthermore, the identity principle for Coleman functions
implies that a Coleman function $f$ vanishes identically if and only if $df=0$
and $f(x)=0$, for some $x\in\{0,1,\infty\}$. This argument will play a crucial
role in the proof of Proposition \ref{relprop1} below.

\segment{furushothm}{}%%%%%%%%%%
We recall Theorem 2.13 of Furusho \cite{FurushoI}.

\subsection*{Theorem (Furusho)}
For $n\geq 2$, $\Li_n^a$ is defined at $z=1$, and the value
  \[
      \zeta(n):=\Li_n^a(1)\in\QQ_p
  \]
  is independent of the choice of $a$. Moreover, $\zeta(n)=0$ if $n$ is even.

\bigskip
\noindent
By Remark \ref{relrem1} and the independence result above, our final result does not depend on the choice of the branch
of the logarithm. From now on, we will therefore ignore the choice of $a$ and
omit it from the notation.

\segment{relprop1}{}%%%%%%%%%%%%%%%
The following proposition states several functional equations for $\Li_2$ and
$\Li_3$ whose complex counterparts correspond to formulas
(1.11), (1.12), (6.10) and (6.4) of Lewin \cite{Lewin}.

\subsection*{Proposition}
The following relations between Coleman functions hold. 
  \begin{gather*}
       \tag{1}\label{eq1}
    \Li_2(z)+\Li_2(1-z)+\log(z)\log(1-z) = 0,\\
       \tag{2}\label{eq2}
    \Li_2(z)+\Li_2 \left( \frac{z}{z-1} \right)+\frac{1}{2}\log(1-z)^2=0,\\
       \tag{3}\label{eq3}
    \Li_3(z)+\Li_3(1-z)+\Li_3 \left( \frac{z}{z-1} \right)-\frac{1}{6}\log(1-z)^3
      +\frac{1}{2}\log(z)\log(1-z)^2 = \zeta(3),\\
       \tag{4}\label{eq4}
    \Li_3(z^2) -4\Li_3(z)-4\Li_3(-z) = 0.
  \end{gather*}
  
 \begin{proof}
 To prove \eqref{eq1} we set
\[
   f(z):=\Li_2(z)+\Li_2(1-z)+\log(z)\log(1-z).
\]
This is a sum of three Coleman function on $Y$. Each of these functions is
defined at $z=0$ and takes the value $0$. It follows that $f(z)$ is
itself a Coleman function, defined at $z=0$, and that $f(0)=0$. 

We compute its derivative with respect to
$z$, using \ref{Lidiffeq}($*$):
\[
    \frac{d}{dz}f(z)=-\frac{\log(1-z)}{z} + \frac{\log(z)}{1-z}
      + \frac{\log(1-z)}{z} -\frac{\log(z)}{1-z} = 0.
\]
Since $f(z)$ is locally analytic, this means that $f(z)$ is constant on every
residue disk. But then the identity principle for Coleman functions implies
that $f(z)$ is constant everywhere. Since $f(0)=0$, $f$ vanishes
identically. This proves \eqref{eq1}. The proof of \eqref{eq2} is very
similar and left to the reader. 

To prove \eqref{eq3} we set
\[
   g(z):=\Li_3(z)+\Li_3(1-z)+\Li_3 \left( \frac{z}{z-1} \right)-\frac{1}{6}\log(1-z)^3
      +\frac{1}{2}\log(z)\log(1-z)^2.
\]
All terms in this sum are Coleman functions which are defined at
$z=0$. Moreover, all terms vanish at $z=0$ except for $\Li_3(1-z)$ which takes
the value $\zeta(3)$. Therefore, $g$ is a Coleman function defined at $z=0$
such that $g(0)=\zeta(3)$. We compute the derivative of $g$:
\[ \tag{5} \label{eq7}
\begin{split}
  \frac{d}{dz}g(z)&=\frac{1}{z}\Li_2(z)+\frac{1}{z-1}\Li_2(1-z)
     -\frac{1}{z(z-1)}\Li_2 \left( \frac{z}{z-1} \right) \\
        & \quad -\frac{1}{2z(z-1)}\log(1-z)^2+\frac{1}{z-1}\log(z)\log(1-z)
\end{split}
\]
From \eqref{eq1} and \eqref{eq2} we obtain expressions for $\Li_2(1-z)$ and
$\Li_2(z/(z-1))$ in terms of $\Li_2(z)$, $\log(z)$ and $\log(1-z)$. Plugging
these expressions into \eqref{eq7}, we obtain, after a short computation, 
\[ \tag{6} \label{eq8}
  \frac{d}{dz}g(z)=0.
\]
As before, this shows that $g$ is constant, and then the identity
$g(0)=\zeta(3)$ proves \eqref{eq3}.

The proof of \eqref{eq4} is again very similar. In addition to the arguments
used already, one has to use that pullbacks of Coleman functions via rational
functions are again Coleman functions. Thus, $\Li_3(z^2)$ is a Coleman
function on $\PP^1(\CC_p)-\{0,\pm 1,\infty\}$ which vanishes for $z=0$.
 \end{proof}
 
 \segment{}{Corollary}%%%%%%%%%%%%%%
 We have 
  \[
      \Li_3(1/2)=\frac{7}{8}\zeta(3)+\frac{1}{6}\log(2)^3.
  \]

\begin{proof}
Evaluation of \ref{relprop1}\eqref{eq3} at $z=1/2$ yields the identity
\[ \tag{$*$} \label{eq9}
    2\Li_3(1/2)+\Li_2(-1)-\frac{1}{3}\log(2) = \zeta(3).
\]
Evaluation of \eqref{eq4} at $z=1$ yields
\[ \tag{$**$} \label{eq10}
  \Li_3(-1) = -\frac{3}{4}.
\]
Combining \eqref{eq9} and \eqref{eq10} gives the desired result.
\end{proof}

%%%%%%%%%%%%%%%%%%%%
\section{$\widetilde{7/8} = 7/8$ via complex polylogarithms}%%%%
%%%%%%%%%%%%%%%%

\segment{}{}%%%%%%  
We now give a proof of the relation
\[
   \Li^U_3(1/2)=\frac{7}{8}\zeta^U(3)+\frac{1}{6}\log^U(2)^3,
\]
given to us by one of the referees, which makes use of Brown's motivic polylogarithms and their complex periods. Let $X = \thrpl$. Given a rational number $z$ with $0<z<1$, we let $_z ch_0$ denote the straight light path in $X(\CC)$ from the tangent vector $\partial/\partial t$ at $0$ to $z$. We also let $dch$ denote the straight line path from $\partial/\partial t$ at $0$ to $-\partial/\partial t$ at $1$. Using the notation of Segment \ref{21b}, we define 
\[
\log^{\mM}(z) := \int_{_z ch_0} 0  = \left[ \Uu{_z P_0}(X), {_z ch_0}, \int_0^z 0 \right]^{\m{B}, \dR},
\]
\[
\Li_n^\mM(z) := \int_{_z ch_0} 0^{n-1}1 = \left[ \Uu{_z P_0}(X), {_z ch_0}, \int_0^z 0^{n-1}1 \right]^{\m{B}, \dR},
\]
and
\[
\ze^\mM(n) := \int_{dch} 0^{n-1}1 = \left[ \Uu{_z P_0}(X), dch, \int_0^z 0^{n-1}1 \right]^{\m{B}, \dR}.
\]

\segment{r28b}{}%%%%%%
As mentioned above, our definition differs slightly from Brown's, who considers the rings of functions
\[
\Oo({_z P_0}(X))
\]
in place of the universal enveloping bimodules, and paths from de Rham to Betti realization in place of our paths from Betti to de Rham realization. Nevertheless, the formalism is quite similar. The complex period map (\ref{19e}) sends our motivic polylogarithms and zeta values to the corresponding complex polylogarithms and zeta values. The component
\[
d: ({_\dR A _\m{B}})_n \to ( A_\dR )_1 \otimes  ({_\dR A _\m{B}})_{n-1}
\]
of the coaction satisfies the formulas
\[
d \log^\mM(z) = \log^U(z) \otimes 1,
\]
\[
d\Li^\mM_n(z) = \log^U(z) \otimes \Li^\mM_{n-1}(z),
\]
and
\[
\Li^\mM_1(z) = -\log^\mM(1-z).
\]
Moreover, if $\xi \in ({_\dR A _\m{B}})_n$ is homogeneous of degree $n >1$ and $d\xi=0$, then $\xi \in \QQ \ze^\mM(n)$ (c.f. \cite[Theorem 2.7]{BrownICM}).

\segment{r26a}{Lemma}%%%%%%%%%%%
We have
\[
\Li ^\mM_3(\frac{1}{2}) = \frac{7}{8} \ze^\mM(3) + \frac{1}{6}\log^\mM(2)^3  - \frac{1}{2} \log^\mM(2) \ze^\mM(2).
\]

\begin{proof}
By the formulas of segment \ref{r28b}, we have
\[
d\Li^\mM_2(1/2) = -\log^U(2) \otimes \log^\mM(2)
\]
and
\[
d\log^\mM(2)^2 = 2\log^U(2) \otimes \log^\mM(2).
\]
So we have
\[
\Li_2^\mM(1/2) = -\frac{1}{2} \log^\mM(2)^2 + \al \ze^\mM(2)
\]
for some $\al \in \QQ$. Applying the period map and comparing with the classical identity of complex polylogarithms
\[
\Li^\infty_2(1/2) = -\frac{1}{2} \log^\infty(2)^2 + \ze^\infty(2)/2,
\]
we obtain
\[
\al = \frac{1}{2}.
\]

Moving up to level $3$, we compute
\begin{align*}
d\Li_3^\mM(1/2) 	&= -\log^U(2) \otimes \Li^\mM_2(1/2) \\
		&= \frac{1}{2} \log^U(2) \otimes \log^\mM(2)^2
		-\frac{1}{2} \log^U(2) \otimes \ze^\mM(2) \\
		&= d \left( \frac{1}{6} \log^\mM(2)^3 -\frac{1}{2} \log^\mM(2) \ze^\mM(2)  \right),			
\end{align*}
from which
\[
\Li_3^\mM(1/2) = \al' \ze^\mM(3) + \frac{1}{6}\log^\mM(2)^3 - \frac{1}{2} \log ^\mM(2) \ze^\mM(2)
\]
for some $\al' \in \QQ$. We now apply the period map again and compare with the identity of complex polylogarithms  \ref{r28b}($*$),  to obtain
\[
\al' = 7/8.
\]
\end{proof}

Pulling back along the orbit map $o({_\dR 1_\m{B}^+})$ of Segment \ref{21b}, we obtain the identity of unipotent motivic polylogarithms
\[
\Li_3^U(1/2) = \frac{7}{8}\ze^U(3) + \frac{1}{6} \log ^U(2)^3.
\]
As mentioned above, we also obtain the identity of $p$-adic polylogarithms for every prime $p$ by applying the $p$-adic period map.

%qqq

\bibliography{references}

\bibliographystyle{alphanum}

\vfill

\Small\textsc{I.D. Fakult\"at f\"ur Mathematik, Universit\"at Duisburg-Essen, Thea-Leymann-Strasse 9, 45127 Essen, Germany}, {E-mail address:} \texttt{ishaidc@gmail.com}

\Small\textsc{S.W. Institut f\"ur Reine Mathematik, Universit\"at Ulm, Helmholtzstrasse 18, 89081 Ulm, Germany}

\end{document}